\documentclass[aps,pre,reprint,superscriptaddress,longbibliography]{revtex4-2}
\usepackage[T1]{fontenc}
\usepackage{graphicx}
\usepackage{amsmath,amssymb,amsfonts}

\usepackage{amsthm}
\usepackage{mathrsfs}
\usepackage{booktabs}
\usepackage[hidelinks]{hyperref}
\theoremstyle{plain}
\newtheorem{theorem}{Theorem}
\newtheorem{proposition}[theorem]{Proposition}
\newtheorem{corollary}[theorem]{Corollary}
\theoremstyle{remark}
\newtheorem{remark}{Remark}
\theoremstyle{definition}
\newtheorem{definition}{Definition}
\newcommand{\bfone}{\mathbf{1}}
\newcommand{\cump}{\operatorname{cum}}
\newcommand{\diag}{\operatorname{diag}}
\newcommand{\tr}{\operatorname{tr}}
\newcommand{\Etwo}{\mathcal E_N}
\newcommand{\dd}{\,{\rm d}}
\newcommand{\bE}{\mathbb E}
\newcommand{\bP}{\mathbb P}
\newcommand{\bR}{\mathbb R}
\newcommand{\bT}{\mathbb T}
\newcommand{\Var}{\operatorname{Var}}
\newcommand{\Cov}{\operatorname{Cov}}
\newcommand{\cF}{\mathcal F}
\newcommand{\cL}{\mathcal L}
\newenvironment{Thm}[1][]{\begin{theorem}[#1]}{\end{theorem}}
\newenvironment{Prop}[1][]{\begin{proposition}[#1]}{\end{proposition}}
\newenvironment{Cor}[1][]{\begin{corollary}[#1]}{\end{corollary}}
\newenvironment{Def}[1][]{\begin{definition}[#1]}{\end{definition}}
\newenvironment{Rem}[1][]{\begin{remark}[#1]}{\end{remark}}

\begin{document}
\title{An exchangeable pair-Hawkes reference model and a cumulant obstruction
for hard-sphere collision statistics}
\author{Luis Iv\'an Hern\'andez Ru\'iz}
\email{luisivanhr@ciencias.unam.mx}
\affiliation{AI Division, AlpacaTech Co., Ltd., H1O Hirakawacho 703,
1-6-4 Hirakawacho, Chiyoda-ku, Tokyo 102-0093, Japan}
\begin{abstract}
Collision counts in a dilute gas of hard particles reflect the dependence
between successive encounters. We study whether these fluctuations can be
represented by a positive linear Hawkes process, in which each event adds a
nonnegative contribution to the rate of future events. We construct such a
model for particle pairs, require invariance under particle relabeling, and
derive its stationary mean and covariance. We prove that, for a stationary
positive Hawkes process with finitely many event types and finite expected
total family sizes, the long-time growth rate of the third centered moment of
any count formed by selecting event types is at least that of its variance.
More generally, the long-time growth rates are nondecreasing with cumulant
order. These predicted growth rates therefore cannot be matched
simultaneously to the ones found in literature in approximations that follow
velocity and collision count for a particle in an equilibrium gas.
Simulations support the covariance predictions and illustrate the dependence
of the physical comparison on the observation window. This restriction
motivates models that retain the state left by a collision or allow past
events to reduce future rates.
\end{abstract}
\maketitle

\section{Introduction}\label{sec:introduction}

A collision record specifies the times of contact and the pairs involved.
Each collision changes the particles' velocities and hence the waiting
times to later encounters. We ask how much of this dependence can be
described by a model that records event times and particle labels. The mean
collision rate gives the average activity. The variance and higher
cumulants describe its fluctuations; for a count, the first cumulant is
the mean, the second is the variance, and the third is the third centered
moment. A Poisson reference process gives equal values for these
quantities. For a single particle followed through its collisions, called
a tagged particle, the dilute-gas calculation of Visco, van Wijland and
Trizac \cite{ViscoVanWijlandTrizac2008} gives different values. The
relations between the cumulants can therefore help distinguish models
with the same mean collision rate.

We consider a positive linear Hawkes process, in which every past event
adds a nonnegative contribution to the present event rate. This provides
a simple description of event clustering. To test its range of possible
fluctuations, we compare the long-time growth rates of its cumulants. If
\(N_A(0,t]\) counts the events whose types belong to a nonempty subset
\(A\), write
\[
 c_j^A=\lim_{t\to\infty}t^{-1}\cump_j\{N_A(0,t]\}.
\]
For pair-labeled events, choosing all pairs that contain a given particle
produces its tagged count. We prove that a stationary subcritical positive
Hawkes process with finitely many types satisfies
\(c_{j+1}^A\geq c_j^A\) for every \(j\geq1\). Subcriticality means
that the family generated by an event has finite expected total size.
The cluster representation explains the inequality: independently arriving
immigrants start families of descendants, and each family contributes a
nonnegative integer number of events to the selected count. Successive
positive integer powers of this number are nondecreasing. The cumulant
formula transfers that ordering to the long-time rates.

To compare this restriction with gas dynamics, we use a kinetic equation
that follows a tagged particle's velocity together with its collision
count. Visco, van Wijland and Trizac obtain explicit cumulants by
approximating an auxiliary velocity profile that weights trajectories
according to their collision counts. They use a Gaussian profile multiplied
by one plus a polynomial correction in squared speed, allowing its width
and the correction coefficient to vary. This is the first Sonine
approximation, described in
Sec.~\ref{sec:sonine-benchmark}. It reduces the velocity-dependent
equation to relations between a few moments. Its cumulant rates,
denoted by \(c_j^{\rm S}\), satisfy
\(c_3^{\rm S}<c_2^{\rm S}\) for every \(d\geq2\). Our inequality
therefore rules out these approximate kinetic cumulants as the exact
long-time cumulants of a positive Hawkes subset count. The comparison
uses a theorem for the point process and a specified approximation for
the gas.

Hawkes introduced this process in \cite{Hawkes1971Spectra}.
Hawkes and Oakes
\cite{HawkesOakes1974} represented this process as a collection of independent
clusters, and general cumulant formulas were obtained in
\cite{JovanovicHertzRotter2015}. Stability, covariance formulas and limit
theorems are treated in
\cite{BremaudMassoulie1996,BremaudMassoulie2001,
BacryDelattreHoffmannMuzy2013,BacryMuzy2016,Zhu2013}.
Further accounts and extensions may be found in
\cite{DaleyVereJones2003,BacryMastromatteoMuzy2015,
DelattreFournierHoffmann2016,HernandezRuiz2025RenewalHawkesLLN,
HernandezYano2025Cluster}.

For a gas with \(N\) particles, the natural index set for collision
labels is
\[
 \Etwo=\{\{i,j\}:1\leq i<j\leq N\}.
\]
Two pair labels may coincide, share one particle, or be disjoint. We impose
exchangeability, meaning that relabeling the particles leaves the process
unchanged in law. Its interaction kernels then depend only on these three
relations. This symmetry separates fluctuations into total activity,
differences between tagged particle counts, and weighted differences of
pair counts whose contributions cancel when summed around any particle.
We call these the global,
tagged particle and cycle modes. Working separately on these modes gives
the stationarity condition and covariance formulas.

The organization of pairs by these relations is known as the Johnson
scheme \(J(N,2)\). Matrices whose entries are constant on each relation
form its Bose--Mesner algebra
\cite{BoseMesner1959,Delsarte1973,BrouwerCohenNeumaier1989,Bailey2004}.
The shared-particle relation also describes the line graph of the complete
graph: its vertices are particle pairs, with adjacent vertices sharing a
particle \cite{GodsilRoyle2001}. These equivalent descriptions let us use
the classical algebraic decomposition to calculate fluctuations of the
collision record.

Between collisions, hard-spheres retain their velocities and move freely.
Each collision changes the two velocities while preserving total energy and
momentum. Kinetic theory describes this motion under specified dilute-gas
assumptions
\cite{Kac1956,Cercignani1988,ChapmanCowling1970,
CercignaniIllnerPulvirenti1994,Spohn1991}.
At fixed volume, the dilute scaling in which the particle diameter
\(a\) tends to zero and \(N\) grows with \(Na^{d-1}\) tending to a
positive constant is called the Boltzmann--Grad regime. The derivation
of kinetic equations in this regime is discussed in
\cite{Grad1958,Lanford1975,GallagherSaintRaymondTexier2014,
Saffirio2016,PulvirentiSimonella2017}; equilibrium fluctuations require
further estimates
\cite{BodineauGallagherSaintRaymondSimonella2020,
BodineauGallagherSaintRaymondSimonella2024}.
Here we use the pair-Hawkes process as a reference law for event times and
pair labels. Kinetic theory supplies the mean collision rate used for
calibration and the Sonine cumulants used for comparison.

The cumulant restriction also motivates changes to the dependence between
events. A state retained after each event can affect the waiting time to
the next one. We give a jump process with two states whose count satisfies
\(c_3<c_2\), and relate the velocity of a tagged particle to this type of
memory. Another possibility is inhibition: allowing a past event to reduce
the present rate. We examine a signed Hawkes model whose intensity is
clipped at zero to keep it nonnegative, and prove a sufficient condition
for stationarity.

Section~\ref{sec:model} defines the collision counts and the reference
process. Section~\ref{sec:obstruction} proves the cumulant inequality and
compares it with the hard-particle benchmark. The Johnson decomposition,
mean-rate calibration and covariance formulas are given in
Secs.~\ref{sec:johnson}--\ref{sec:covariance}. Numerical comparisons
with the covariance formulas and hard-particle cumulants are presented
in Sec.~\ref{sec:numerical}.
Section~\ref{sec:extensions} considers the extensions, and
Sec.~\ref{sec:conclusion} gives the conclusions.

\section{Collision counts and the pair-Hawkes process}\label{sec:model}

\subsection{Hard sphere dynamics and its counting measure}

Let \(d\geq2\) and \(N\geq4\). Consider \(N\) identical hard-spheres of
diameter \(0<a<L/2\) and mass \(m>0\) in the periodic box \(\bT_L^d\). The
admissible
phase space is
\[
 \Omega_N=\left\{(q,v)\in(\bT_L^d\times\bR^d)^N:
 |q_i-q_j|_{\bT_L^d}\geq a,\ i\neq j\right\}.
\]
Outside contact times, \(\dot q_i=v_i\) and \(\dot v_i=0\). At a regular
binary contact of particles \(i\) and \(j\), let \(r_{ij}\) be the unique
minimal image representative of \(q_i-q_j\) in
\([-L/2,L/2)^d\). Then \(|r_{ij}|=a\); set \(n=r_{ij}/a\) and assume that
\((v_i^--v_j^-)\cdot n<0\). The postcollisional velocities are
\begin{equation}\label{eq:collision-rule}
 v_i^+=v_i^- -\{(v_i^--v_j^-)\cdot n\}n,\qquad
 v_j^+=v_j^- +\{(v_i^--v_j^-)\cdot n\}n .
\end{equation}
This transformation preserves \(\sum_i mv_i\) and
\(\sum_i m|v_i|^2/2\). Simultaneous multiple contacts form a singular set and
are excluded from the collision records considered below.

For \(e=\{i,j\}\in\Etwo\), let \(\{T_{e,k}\}_{k\in\mathbb Z}\) be the locally
finite set of regular collision times of the pair \(e\). Its random counting
measure is
\begin{equation}\label{eq:counting-measure}
 N_e(A)=\sum_k\mathbf 1_{\{T_{e,k}\in A\}},\qquad
 N_e(\dd s)=\sum_k\delta_{T_{e,k}}(\dd s)
\end{equation}
for Borel sets \(A\subset\bR\). In particular,
\[
 \int_{(-\infty,t)}g(t-s)N_e(\dd s)
 =\sum_{k:T_{e,k}<t}g(t-T_{e,k}).
\]
We use \(N_e(\dd s)\), or equivalently \(\dd N_e(s)\), for the random counting
measure associated with the collision times of the pair \(e\). The total and
tagged counts in a time interval \(I\) are
\begin{equation}\label{eq:physical-counts}
 C_{\rm tot}(I)=\sum_{e\in\Etwo}N_e(I),\qquad
 C_i(I)=\sum_{j\neq i}N_{\{i,j\}}(I).
\end{equation}

\subsection{Conditional intensities and exchangeability}

Let \(\cF_t\) be the history generated by all pair counts up to time \(t\).
A predictable, nonnegative process \(\lambda_e(t)\) is an
\(\cF_t\)-conditional intensity of \(N_e\) if
\[
 M_e(t)=N_e((0,t])-\int_0^t\lambda_e(s)\dd s
\]
is a local martingale. If \(\lambda_e\) has almost surely locally bounded
left continuous paths, the corresponding small time characterization is
\[
 \bP\{N_e((t,t+\Delta])=1\mid\cF_{t^-}\}
   =\lambda_e(t)\Delta+o(\Delta),
 \qquad \Delta\downarrow0,
\]
while the conditional probability of two or more events is \(o(\Delta)\);
see \citep{Bremaud1981,DaleyVereJones2003}.

\begin{Def}[Exchangeable pair-Hawkes reference process]\label{def:pair-hawkes}
Let \(h_{\rm same},h_{\rm share},h_{\rm dis}\) be nonnegative integrable
functions on \([0,\infty)\), and let \(\mu>0\). The exchangeable
pair-Hawkes process is an \(\Etwo\)-variate point process whose ground process
is simple, namely
\[
 \sum_{e\in\Etwo}N_e(\{t\})\leq1
 \quad\hbox{for every \(t\), almost surely},
\]
and whose conditional intensities are
\begin{align}
\lambda_{ij}(t)={}&\mu+
 \int_{(-\infty,t)}h_{\rm same}(t-s)N_{ij}(\dd s)\nonumber\\
&+\sum_{k\neq i,j}\int_{(-\infty,t)}h_{\rm share}(t-s)\nonumber\\
&\hspace{17mm}\times\{N_{ik}(\dd s)+N_{jk}(\dd s)\}\nonumber\\
&+\sum_{\substack{k<\ell\\\{k,\ell\}\cap\{i,j\}=\varnothing}}
 \int_{(-\infty,t)}h_{\rm dis}(t-s)N_{k\ell}(\dd s).
\label{eq:pair-hawkes}
\end{align}
\end{Def}

The kernel \(h_{\rm same}\) describes how an event affects later events of
the same pair. The kernel \(h_{\rm share}\) carries this influence through a
shared particle. A simple choice for a dilute reference model is
\(h_{\rm dis}=0\). Even then, disjoint pairs can become dependent through a
chain of pairs with shared particles. A nonzero \(h_{\rm dis}\) can describe
a common density or activity fluctuation. Its value is a choice of the
reference model.

For a permutation \(\pi\) of \(\{1,\ldots,N\}\), write
\(\pi\{i,j\}=\{\pi(i),\pi(j)\}\). A pair process is exchangeable when
\(\{N_{\pi e}(A):e,A\}\) has the same law as
\(\{N_e(A):e,A\}\). Equation \eqref{eq:pair-hawkes}, together with a common
baseline and a unique stationary law, is invariant under this action.

\begin{Def}[Stationarity]\label{def:stationarity}
A pair point process on \(\bR\) is stationary if, for every \(u\in\bR\),
\(\{N_e(A+u):e,A\}\) and \(\{N_e(A):e,A\}\) have the same finite dimensional
distributions.
\end{Def}

\begin{Rem}[Stationarity and physical equilibrium]
\label{rem:stationarity-equilibrium}
Stationarity means that shifting the time origin leaves the law of the
record unchanged \citep{Kallenberg2017}. Both an equilibrium gas and a driven
nonequilibrium steady state can have this property
\citep{EvansMorriss2008}. Thermodynamic equilibrium imposes further
conditions on the physical system.
\end{Rem}

The state recorded by Definition \ref{def:pair-hawkes} consists of pair
labels and event times. The positions and velocities in
\eqref{eq:collision-rule} belong to the underlying mechanical description.

\section{The cumulant obstruction and the hard-particle benchmark}
\label{sec:obstruction}
\subsection{Nonnegative subset counts}

Let \(N=(N_1,\ldots,N_m)\) be a stationary positive multitype Hawkes process
with immigrant rates \(\mu_1,\ldots,\mu_m\), integrated offspring matrix
\(K\), and \(\rho(K)<1\). In the Poisson cluster construction of
\citet{HawkesOakes1974}, the immigrants form independent Poisson processes on
\(\bR\) \citep{Kingman1993}, and each immigrant carries an independent
subcritical multitype Galton--Watson cluster \citep{AthreyaNey1972}.

For a nonempty subset \(A\subset\{1,\ldots,m\}\), define
\[
 N_A(0,t]=\sum_{e\in A}N_e((0,t]).
\]
Let \(S_A^{(f)}\) denote the total number of points with type in \(A\) in a
complete cluster whose immigrant has type \(f\), including the immigrant when
\(f\in A\).

\begin{Thm}[Subset count cumulant formula and obstruction]\label{thm:cumulant}
Fix \(n\geq1\) and assume
\(\bE[(S_A^{(f)})^{n+1}]<\infty\) for every \(f\). Then the long time
cumulant rates of orders \(j=n,n+1\) exist and satisfy
\begin{equation}\label{eq:cluster-cumulant-rate}
\begin{split}
 c_j^A&:=\lim_{t\to\infty}\frac1t\cump_j\{N_A(0,t]\}\\
 &=\sum_{f=1}^m\mu_f\bE[(S_A^{(f)})^j],
 \qquad j\in\{n,n+1\}.
\end{split}
\end{equation}
Consequently
\begin{equation}\label{eq:cumulant-monotonicity}
 c_{n+1}^A\geq c_n^A.
\end{equation}
In particular, \(c_3^A\geq c_2^A\).
\end{Thm}

\begin{proof}
Consider a cluster with immigrant type \(f\), relative event times
\(\tau_j\geq0\), and event types \(E_j\). If its immigrant is placed at
\(u\in\bR\), its contribution to the window is
\[
 Y_{A,t}^{(f)}(u)=
 \sum_j\mathbf 1_{\{E_j\in A\}}
\mathbf 1_{\{0<u+\tau_j\leq t\}}.
\]
The assumed finite \((n+1)\)-st moment ensures that the cluster contains
finitely many points of the selected types almost surely. We add the
contributions of all immigrant clusters using the exponential formula for
marked Poisson processes:
\begin{equation}\label{eq:poisson-log-characteristic}
\begin{split}
 &\log\bE\exp\{i\theta N_A(0,t]\}\\
 &\quad=\sum_{f=1}^m\mu_f\int_{\bR}
 \left[\bE\exp\{i\theta Y_{A,t}^{(f)}(u)\}-1\right]\dd u.
\end{split}
\end{equation}
See \citep{Kingman1993,LastPenrose2017,Kallenberg2017} for the Poisson random
measure formula. To differentiate it, we first bound the contribution of a
fixed finite cluster. Expanding \((Y_{A,t}^{(f)}(u))^j\) produces at most
\((S_A^{(f)})^j\) ordered tuples. A tuple contributes precisely when all its
points fall inside the window. The possible immigrant locations therefore
form an interval of length at most \(t\), and
\[
 \int_{\bR}(Y_{A,t}^{(f)}(u))^j\dd u
 \leq t(S_A^{(f)})^j.
\]
The moment assumption makes the expectation of the right side finite
for \(j=n,n+1\). Moreover, for real \(\theta\),
\[
 \left|\frac{\partial^j}{\partial\theta^j}
 e^{i\theta Y_{A,t}^{(f)}(u)}\right|
 =(Y_{A,t}^{(f)}(u))^j.
\]
Fubini's theorem and dominated differentiation may therefore be applied to
\eqref{eq:poisson-log-characteristic}. Differentiating \(j\) times at
\(\theta=0\) and multiplying by \(i^{-j}\) gives
\begin{equation}\label{eq:poisson-cumulant}
 \cump_j\{N_A(0,t]\}
 =\sum_{f=1}^m\mu_f\int_{\bR}
 \bE[(Y_{A,t}^{(f)}(u))^j]\dd u.
\end{equation}
Integration over all immigrant locations in \(\bR\) includes clusters that
cross either endpoint of the observation window. Equation
\eqref{eq:poisson-cumulant} therefore keeps their contributions.

For one fixed finite cluster, consider an ordered \(j\)-tuple of
\(A\)-points. All these points lie in the window exactly when the immigrant
location belongs to the intersection of the intervals
\((-\tau_{j_q},t-\tau_{j_q}]\). Summing the lengths of these intersections
gives
\begin{equation}\label{eq:cluster-window-identity}
\begin{split}
 &\frac1t\int_{\bR}(Y_{A,t}^{(f)}(u))^j\dd u\\
 &\quad=\sum_{\substack{\ell_1,\ldots,\ell_j\\E_{\ell_q}\in A}}
 \left(1-\frac{\max_q\tau_{\ell_q}-\min_q\tau_{\ell_q}}t\right)_+.
\end{split}
\end{equation}
Every summand converges to one, while the right-hand side is bounded by
\((S_A^{(f)})^j\). Dominated convergence, first with \(j=n\) and then with
\(j=n+1\), proves both identities in
\eqref{eq:cluster-cumulant-rate}. Finally, \(S_A^{(f)}\) is a nonnegative
integer and therefore
\((S_A^{(f)})^{n+1}\geq(S_A^{(f)})^n\) pointwise. Multiplication by
\(\mu_f\) and summation over \(f\) give
\eqref{eq:cumulant-monotonicity}.
\end{proof}

\begin{Rem}[Moment assumption and scope]\label{rem:scope}
A finite-type Hawkes process with Poisson offspring and \(\rho(K)<1\)
satisfies the moment assumption. To verify this, consider the probability
generating functions \(F_f(z)=\bE[z^{S_A^{(f)}}]\). The immigrant is counted
when its type belongs to \(A\), and its
offspring generate independent subclusters. The branching property gives
the finite system
\[
\begin{gathered}
 F_f(z)=z^{\mathbf 1_{\{f\in A\}}}
 \exp\left\{\sum_{e=1}^mK_{ef}\{F_e(z)-1\}\right\},\\
 1\leq f\leq m.
\end{gathered}
\]
At \(z=1\), one has \(F_f(1)=1\), and the Jacobian of the right side
with respect to \(F=(F_1,\ldots,F_m)\) is \(K^\top\). Since
\(\rho(K^\top)<1\), \(I_m-K^\top\) is invertible. More explicitly, if
\(\mathcal T(z,F)\) denotes the right side and
\(\mathcal H(z,F)=F-\mathcal T(z,F)\), then
\[
 D_F\mathcal H(1,\bfone_m)=I_m-K^\top.
\]
The analytic implicit function theorem therefore gives a unique analytic
solution near
\((z,F)=(1,\bfone_m)\). This branch is the probability generating function.
Indeed, the subcritical total progeny is finite almost surely. The actual
generating function is therefore defined for \(0\leq z\leq1\), satisfies the
displayed system, and converges componentwise to \(\bfone_m\) as
\(z\uparrow1\). For \(z<1\) sufficiently close to one, local uniqueness
identifies it with the analytic branch. Its left derivatives at one are
finite and, by monotone convergence applied to
\(\bE[(S_A^{(f)})_kz^{S_A^{(f)}-k}]\), give every factorial moment. Thus, all
ordinary moments are finite; see \citep{Harris1963,AthreyaNey1972}.
Total and tagged particle counts add events with unit weight, so Theorem
\ref{thm:cumulant} applies to them. A cycle observable uses both positive and
negative coefficients; its cluster contribution can have either sign.
\end{Rem}

\subsection{A kinetic approximation for the tagged collision count}
\label{sec:sonine-benchmark}

The Hawkes inequality can be tested against a kinetic description that
retains the velocity left by each collision. Visco, van Wijland and Trizac
\cite{ViscoVanWijlandTrizac2008} considered the number
\(\mathcal N(t)\) of collisions of a tagged particle in a dilute
equilibrium gas of hard particles. Their collision-counting Boltzmann
equation follows the joint dependence on velocity and collision number.
Multiplying the contribution of a trajectory by
\(e^{s\mathcal N(t)}\), with counting parameter \(s\), gives an auxiliary
velocity profile. Its total mass is the moment generating function of the
collision count, so its long-time exponential growth determines the
cumulant rates.

To obtain explicit rates, they approximate the normalized long-time
profile by a Gaussian multiplied by one plus a Sonine polynomial
correction in squared speed. Sonine polynomials are orthogonal for the
weight induced by the Gaussian radial distribution. They organize the
corrections through velocity moments. The Gaussian width accounts for the
second moment, and the leading nonzero correction retained here is
controlled by the fourth moment. This is the first Sonine approximation.
Taking the zeroth, second and fourth velocity moments of the kinetic
equation gives a small system for the Gaussian width, the correction
coefficient and the growth rate. Expanding its solution near \(s=0\) gives
the mean rate, variance rate and third cumulant rate. The profile being
approximated belongs to the collision-count generating function; the
equilibrium velocity distribution in this kinetic description is
Maxwellian.

Write \(c_j^{\rm S}\) for the long-time cumulant rates obtained by this
approximation. Their Eqs.~(50a)--(50c) give
\begin{align}
\frac{c_1^{\rm S}}{\omega}={}&1,\label{eq:sonine1}\\
\frac{c_2^{\rm S}}{\omega}={}&
\frac9{64}\left(8+\frac1{4d+3}\right),\label{eq:sonine2}\\
\frac{c_3^{\rm S}}{\omega}={}&
\frac{28d\{64d(320d+729)+35775\}+257391}
{8192(4d+3)^3}.\label{eq:sonine3}
\end{align}
Here \(\omega\) is the mean tagged particle collision frequency.
Sonine polynomial approximations of this type are standard in kinetic
treatments of dilute and granular gases
\citep{ChapmanCowling1970,BrilliantovPoschel2004}. A direct
algebraic subtraction gives
\begin{equation}\label{eq:sonine-gap}
\begin{split}
\frac{c_3^{\rm S}-c_2^{\rm S}}{c_1^{\rm S}}
&=-\frac{16384d^3+39168d^2+21276d+1809}
{8192(4d+3)^3}\\
&<0,\qquad d\geq2.
\end{split}
\end{equation}
Every term in the numerator is positive for \(d\geq2\). The displayed
quantity is therefore negative throughout these dimensions: the first Sonine
benchmark satisfies \(c_3^{\rm S}<c_2^{\rm S}\).

Theorem~\ref{thm:cumulant} requires the opposite ordering for every
nonnegative subset count of a stationary subcritical positive Hawkes process
with finitely many event types. The first Sonine coefficients in
Eqs.~\eqref{eq:sonine1}--\eqref{eq:sonine3} are therefore incompatible with
the exact long-time cumulants of such a count. This includes the tagged
count in Definition~\ref{def:pair-hawkes}. The corresponding sign for exact
hard-sphere dynamics remains a separate question.

For comparison, a scalar Hawkes process with branching ratio \(\eta\) has
\[
 \frac{c_2}{c_1}=\frac1{(1-\eta)^2},\qquad
 \frac{c_3}{c_1}=\frac{1+2\eta}{(1-\eta)^4}.
\]
Matching the second Sonine cumulant in dimension three fixes
\(\eta=1-\sqrt{320/363}\). With this choice, the third Hawkes cumulant exceeds
the Sonine value. The subset theorem shows that the discrepancy in ordering
persists for any finite number of positive event types.

The numerical comparisons in Sec.~\ref{sec:numerical} place this
benchmark alongside finite-window collision statistics. The planar results
support the reversed ordering at the longer observation window, and the
three-dimensional calculation leaves the ordering unresolved.

\section{The Johnson scheme reduction}\label{sec:johnson}

We next apply exchangeability to the pair process. Set
\(M=\binom N2\), and define the kernel masses by
\[
 k_\alpha=\int_0^\infty h_\alpha(t)\dd t,\qquad
 \alpha\in\{{\rm same},{\rm share},{\rm dis}\},
\]
Let \(K=(K_{ef})_{e,f\in\Etwo}\) denote the integrated kernel matrix. Hence
\begin{equation}\label{eq:integrated-kernel}
K_{ef}=
\begin{cases}
k_{\rm same},&e=f,\\
k_{\rm share},&|e\cap f|=1,\\
k_{\rm dis},&e\cap f=\varnothing.
\end{cases}
\end{equation}

\begin{Prop}[Exchangeability and the Bose--Mesner algebra]\label{prop:BM}
A matrix \(A=(A_{ef})_{e,f\in\Etwo}\) satisfies
\(A_{\pi e,\pi f}=A_{ef}\) for every particle permutation \(\pi\) if and only
if it is constant on the three relations \(e=f\), \(|e\cap f|=1\), and
\(e\cap f=\varnothing\). Consequently, \(K\) in
\eqref{eq:integrated-kernel} belongs to the Bose--Mesner algebra of
\(J(N,2)\).
\end{Prop}

\begin{figure}[htbp]
\centering
\includegraphics[width=\columnwidth]{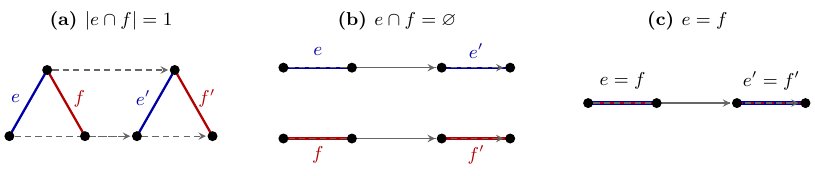}
\caption{The endpoint map used in Proposition~\ref{prop:BM}. The dashed
arrows define \(\varphi:e\cup f\to e'\cup f'\) in the three possible cases.
After \(\varphi\) has been constructed, the unused particle labels can be
matched by any bijection.}
\label{fig:endpoint-bijection}
\end{figure}

\begin{proof}
The cardinality \(|e\cap f|\in\{0,1,2\}\) is preserved by every permutation.
Conversely, suppose that \(|e\cap f|=|e'\cap f'|\). If the intersection has
size one, define \(\varphi:e\cup f\to e'\cup f'\) by sending the common
vertex to the common vertex and the remaining endpoints edge by edge. If the
intersection has size zero, map the two endpoints of \(e\) bijectively to
those of \(e'\) and the two endpoints of \(f\) bijectively to those of
\(f'\). If it has size two, then \(e=f\) and \(e'=f'\), and any bijection
between their two endpoints can be used. These three endpoint maps are shown
in Fig.~\ref{fig:endpoint-bijection}. In each case, distinct used vertices
have distinct images and satisfy
\[
 |e\cup f|=4-|e\cap f|=4-|e'\cap f'|=|e'\cup f'|.
\]
Thus, \(\varphi\) is a bijection. The complements of its domain and range
have the same cardinality, so choose a bijection
\(\psi:[N]\setminus(e\cup f)\to[N]\setminus(e'\cup f')\). Then
\(\pi=\varphi\cup\psi\) is a bijection from \([N]\) to itself, and hence a
permutation, with \(\pi e=e'\) and \(\pi f=f'\). Thus, the diagonal action
has exactly the three stated orbits. The orbit indicator matrices span the
Bose--Mesner algebra of \(J(N,2)\); see
\citep{BoseMesner1959,Delsarte1973,Bailey2004}.
\end{proof}

Let \(B\) be the \(N\times M\) unoriented incidence matrix,
\[
 B_{i,e}=\mathbf 1_{\{i\in e\}}.
\]
The three mutually orthogonal subspaces are
\begin{equation}\label{eq:three-spaces}
\begin{split}
 U_G&=\operatorname{span}\{\bfone_M\},\\
 U_P&=\{B^\top u:u^\top\bfone_N=0\},\\
 U_C&=\ker B.
\end{split}
\end{equation}
We call these the global, tagged particle and cycle spaces. The global
space describes changes common to all pairs. Incidence vectors and their
centered differences belong to \(U_G\oplus U_P\), so these two spaces
describe the tagged counts in \eqref{eq:physical-counts}. A vector in
\(U_C\) cancels when we sum its pair coordinates around any particle.

\begin{Rem}[Interpretation of the cycle space]\label{rem:cycle-meaning}
Let \(i,j,k,\ell\) be distinct particle labels and write
\(e_{\{r,s\}}\) for the coordinate vector of the pair \(\{r,s\}\). Set
\(N(I):=(N_e(I))_{e\in\Etwo}\). Then
\[
\begin{split}
 q_{i,j,k,\ell}
 &:=e_{\{i,j\}}-e_{\{j,k\}}+e_{\{k,\ell\}}-e_{\{i,\ell\}}\in U_C,\\
 q_{i,j,k,\ell}^{\top}N(I)
 &=N_{ij}(I)-N_{jk}(I)+N_{k\ell}(I)-N_{i\ell}(I).
\end{split}
\]
Each particle label occurs once with each sign, giving
\(Bq_{i,j,k,\ell}=0\). Adding a multiple of \(q_{i,j,k,\ell}\) to the pair
count vector changes these four coordinates and preserves every tagged
total. The word ``cycle'' names this alternating contrast of coordinates.
Its microscopic interpretation requires information about particle motion
beyond the labels and times in the record.
\end{Rem}

\begin{figure}[htbp]
\centering
\includegraphics[width=\columnwidth]{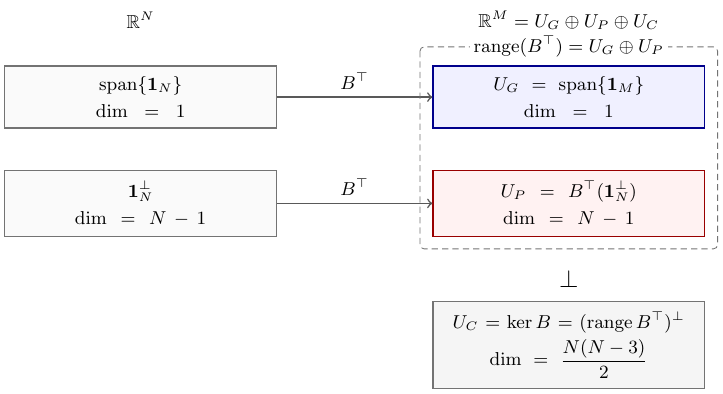}
\caption{Algebraic schematic of the subspaces in
\eqref{eq:three-spaces}. The map \(B^\top\) sends
\(\operatorname{span}\{\bfone_N\}\) and \(\bfone_N^\perp\) onto \(U_G\) and
\(U_P\), respectively, while \(U_C\) is the orthogonal complement of
\(\operatorname{range}(B^\top)\). The direct sum and the dimensions are
proved in Theorem~\ref{thm:three-mode}.}
\label{fig:three-spaces}
\end{figure}

\begin{Thm}[Three mode decomposition]\label{thm:three-mode}
For \(N\geq4\),
\[
 \bR^M=U_G\oplus U_P\oplus U_C,
\]
Moreover
\[
\begin{gathered}
 \dim U_G=1,\qquad \dim U_P=N-1,\\
 \dim U_C=\frac{N(N-3)}2.
\end{gathered}
\]
The matrix \(K\) acts on these spaces by the eigenvalues
\begin{align}
\kappa_G={}&k_{\rm same}+2(N-2)k_{\rm share}
             +\binom{N-2}{2}k_{\rm dis},\label{eq:kappaG}\\
\kappa_P={}&k_{\rm same}+(N-4)k_{\rm share}
             -(N-3)k_{\rm dis},\label{eq:kappaP}\\
\kappa_C={}&k_{\rm same}-2k_{\rm share}+k_{\rm dis}.
\label{eq:kappaC}
\end{align}
\end{Thm}

\begin{proof}
One has
\[
 BB^\top=(N-2)I_N+\bfone_N\bfone_N^\top.
\]
Indeed, the diagonal entries of \(BB^\top\) equal \(N-1\), while two
different rows meet in the unique edge joining their particle labels. The
displayed matrix has eigenvalue \(2(N-1)\) on
\(\operatorname{span}\{\bfone_N\}\) and eigenvalue \(N-2\) on
\(\bfone_N^\perp\). It is therefore nonsingular and
\(\operatorname{rank}B=N\). Moreover
\[
 B^\top\bfone_N=2\bfone_M,\qquad
 B\bfone_M=(N-1)\bfone_N.
\]
The first identity shows that \(U_G\subset\operatorname{range}(B^\top)\).
If \(u^\top\bfone_N=0\), then
\[
 \bfone_M^\top B^\top u=(B\bfone_M)^\top u
 =(N-1)\bfone_N^\top u=0,
\]
so \(U_G\perp U_P\). Since \(B^\top\) is injective,
\(\dim U_P=N-1\), and hence
\[
 \operatorname{range}(B^\top)=U_G\oplus U_P=(\ker B)^\perp.
\]
It follows that \(\dim U_C=M-N=N(N-3)/2\), which proves both the
orthogonal decomposition and the dimension statements.

Let \(A_1\) be the adjacency matrix of the shared particle relation. Since
\((B^\top B)_{ef}=|e\cap f|\), one has
\[
 A_1=B^\top B-2I_M.
\]
The incidence matrix gives its action on each space directly. First,
\[
 A_1\bfone_M
 =B^\top(B\bfone_M)-2\bfone_M
 =2(N-2)\bfone_M.
\]
Second, if \(x=B^\top u\in U_P\), then
\[
 A_1x
 =B^\top(BB^\top u)-2B^\top u
 =\{(N-2)-2\}x=(N-4)x,
\]
because \(u\perp\bfone_N\). Finally, if \(x\in U_C\), then \(Bx=0\)
and \(A_1x=-2x\). Thus, \(A_1\) has the respective eigenvalues
\(2(N-2)\), \(N-4\), and \(-2\) on \(U_G,U_P,U_C\).

Moreover
\[
 K=(k_{\rm same}-k_{\rm dis})I_M
 +(k_{\rm share}-k_{\rm dis})A_1
 +k_{\rm dis}\bfone_M\bfone_M^\top.
\]
The last matrix acts by \(M\) on \(U_G\) and vanishes on
\(U_P\oplus U_C\). Consequently, the three eigenvalues of \(K\) are
\begin{align*}
 & (k_{\rm same}-k_{\rm dis})
 +2(N-2)(k_{\rm share}-k_{\rm dis})+Mk_{\rm dis}\\
 &\quad=k_{\rm same}+2(N-2)k_{\rm share}
 +\binom{N-2}{2}k_{\rm dis},\\
 & (k_{\rm same}-k_{\rm dis})
 +(N-4)(k_{\rm share}-k_{\rm dis})\\
 &\quad=k_{\rm same}+(N-4)k_{\rm share}-(N-3)k_{\rm dis},\\
 & (k_{\rm same}-k_{\rm dis})-2(k_{\rm share}-k_{\rm dis})\\
 &\quad=k_{\rm same}-2k_{\rm share}+k_{\rm dis}.
\end{align*}
These are \eqref{eq:kappaG}--\eqref{eq:kappaC}.
\end{proof}

The orthogonal projectors used below are
\begin{equation}\label{eq:projectors}
\begin{split}
 P_G&=\frac1M\bfone_M\bfone_M^\top,\\
 P_P&=\frac1{N-2}B^\top\left(I_N-\frac1N\bfone_N\bfone_N^\top\right)B,\\
 P_C&=I_M-P_G-P_P.
\end{split}
\end{equation}
To check the formula for \(P_P\), we follow its action on the three spaces.
Put
\(Q_N=I_N-N^{-1}\bfone_N\bfone_N^\top\). Then \(Q_N\bfone_N=0\),
so \(P_P\) annihilates \(U_G\); it also annihilates \(U_C=\ker B\).
For \(B^\top u\in U_P\), with \(u\perp\bfone_N\),
\[
 \frac1{N-2}B^\top Q_NBB^\top u
 =\frac1{N-2}B^\top Q_N\{(N-2)u\}=B^\top u.
\]
The displayed symmetric matrix is therefore the orthogonal projector onto
\(U_P\). The formulas for \(P_G\) and \(P_C\) follow from the orthogonal
decomposition in Theorem \ref{thm:three-mode}.

\begin{Cor}[Stationarity and mean rate]\label{cor:stationarity}
Assume that the kernels in Definition \ref{def:pair-hawkes} are nonnegative.
Then a necessary and sufficient subcriticality condition for the stationary
Poisson cluster construction with finite mean intensity is
\[
 \rho(K)=\kappa_G<1.
\]
In that case the stationary law is unique, nonexplosive and exchangeable, and
every edge has mean intensity,
\begin{equation}\label{eq:mean-rate}
 r_H=\frac{\mu}{1-\kappa_G}.
\end{equation}
\end{Cor}

\begin{proof}
The nonnegative matrix \(K\) has constant row sum \(\kappa_G\). Hence
\(K\bfone_M=\kappa_G\bfone_M\). The maximum row sum bound gives
\(\rho(K)\leq\kappa_G\), while \(\kappa_G\) is itself an eigenvalue; thus
\(\rho(K)=\kappa_G\).

Necessity can be read directly from the first moment equation. If a
stationary solution has finite mean intensity vector
\(r=(r_e)_{e\in\Etwo}\), Tonelli's theorem applied to the nonnegative kernels
in \eqref{eq:pair-hawkes} gives
\[
 r=\mu\bfone_M+Kr.
\]
Since \(K\) is symmetric and has row sum \(\kappa_G\), it also has column sum
\(\kappa_G\). Multiplication by \(\bfone_M^\top\) yields
\[
 (1-\kappa_G)\bfone_M^\top r=M\mu.
\]
The right side is strictly positive, so \(\kappa_G<1\).

Conversely, when \(\kappa_G<1\), the subcritical multitype
Poisson cluster construction gives a unique stationary nonexplosive law with
finite first moments \citep{HawkesOakes1974,DaleyVereJones2003}. Moreover
\[
 r=\mu(I_M-K)^{-1}\bfone_M
   =\mu\sum_{n\geq0}K^n\bfone_M
   =\frac{\mu}{1-\kappa_G}\bfone_M,
\]
which proves \eqref{eq:mean-rate}. Permutation invariance of the specification
and uniqueness in law imply exchangeability.
\end{proof}

\subsection{A finite size parametrization}

Changing \(N\) changes how many pairs share a particle or are disjoint. To
compare these systems, we keep the three modal gains fixed and solve for the
kernel masses. The spectral decomposition is
\[
\begin{split}
 K&=\kappa_GP_G+\kappa_PP_P+\kappa_CP_C\\
 &=\kappa_CI_M+(\kappa_G-\kappa_C)P_G
  +(\kappa_P-\kappa_C)P_P.
\end{split}
\]
For \(e,f\in\Etwo\), the projector entries are
\[
 (P_G)_{ef}=\frac1M,\qquad
 (P_P)_{ef}
 =\frac{|e\cap f|-4/N}{N-2}.
\]
The second identity follows from \(B^\top\bfone_N=2\bfone_M\).
Consequently
\[
 (P_P)_{ef}=
 \begin{cases}
 2/N,&e=f,\\[1mm]
 (N-4)/\{N(N-2)\},&|e\cap f|=1,\\[1mm]
 -4/\{N(N-2)\},&e\cap f=\varnothing.
 \end{cases}
\]
Substitution into the preceding spectral decomposition, relation by
relation, gives
\begin{align}
k_{\rm same}={}&\frac{\kappa_G}{M}+\frac{2\kappa_P}{N}
 +\frac{N-3}{N-1}\kappa_C,\label{eq:inverse-same}\\
k_{\rm share}={}&\frac{\kappa_G-\kappa_C}{M}
 +\frac{(N-4)(\kappa_P-\kappa_C)}{N(N-2)},\label{eq:inverse-share}\\
k_{\rm dis}={}&\frac{\kappa_G-\kappa_C}{M}
 -\frac{4(\kappa_P-\kappa_C)}{N(N-2)}.\label{eq:inverse-dis}
\end{align}
The powers of \(N\) account for the number of pairs in each relation. A fixed
positive mass on disjoint pairs makes \(\kappa_G\) grow as \(N^2\).
At fixed modal gains, \(k_{\rm same}=O(1)\),
\(k_{\rm share}=O(N^{-1})\), and \(k_{\rm dis}=O(N^{-2})\), unless the
corresponding leading coefficient vanishes. A selected triple
\((\kappa_G,\kappa_P,\kappa_C)\) gives a positive Hawkes model only when the
reconstructed entries are nonnegative.

\section{Aggregation and kinetic mean-rate calibration}\label{sec:aggregation}

Summing over all pair labels removes the distinction between the three
relations. Their contributions combine into the global kernel
\begin{equation}\label{eq:hG}
 h_G(t)=h_{\rm same}(t)+2(N-2)h_{\rm share}(t)
       +\binom{N-2}{2}h_{\rm dis}(t).
\end{equation}

\begin{Prop}[Exact scalar aggregation]\label{prop:scalar-aggregation}
For the positive pair-Hawkes process,
\[
 \Lambda_{\rm tot}(t):=\sum_{e\in\Etwo}\lambda_e(t)
 =M\mu+\int_{(-\infty,t)}h_G(t-s)C_{\rm tot}(\dd s).
\]
Consequently, the unmarked total count is a scalar Hawkes process with
baseline \(M\mu\) and kernel \(h_G\).
\end{Prop}

\begin{proof}
Fix a previous edge \(f\). Among the possible target edges, one equals \(f\),
\(2(N-2)\) share one particle with \(f\), and \(\binom{N-2}{2}\) are
disjoint from \(f\). Hence, for every \(u\geq0\),
\[
\begin{split}
 \sum_{e\in\Etwo}h_{ef}(u)
 &=h_{\rm same}(u)+2(N-2)h_{\rm share}(u)\\
 &\quad+\binom{N-2}{2}h_{\rm dis}(u)=h_G(u).
\end{split}
\]
Because \(\Etwo\) is finite, the sum and stochastic integrals can be
interchanged. Summing \eqref{eq:pair-hawkes} first over \(e\) and then over
the source label \(f\) gives
\begin{align*}
&\sum_{e\in\Etwo}\lambda_e(t)\\
&\quad=M\mu+\sum_{f\in\Etwo}\int_{(-\infty,t)}
       \left\{\sum_{e\in\Etwo}h_{ef}(t-s)\right\}N_f(\dd s)\\
&\quad=M\mu+\sum_{f\in\Etwo}\int_{(-\infty,t)}
       h_G(t-s)N_f(\dd s)\\
&\quad=M\mu+\int_{(-\infty,t)}h_G(t-s)C_{\rm tot}(\dd s).
\end{align*}
The ground process is simple by Definition \ref{def:pair-hawkes}, and hence
\(C_{\rm tot}\) is simple. Let \(\mathcal G_t\) be the natural filtration of
\(C_{\rm tot}\). The last member of the display depends only on the past of
\(C_{\rm tot}\), so \(\Lambda_{\rm tot}\) is
\(\mathcal G_{t^-}\)-predictable. The \(\cF_t\)-compensator of
\(C_{\rm tot}\) is
\(\int_0^t\Lambda_{\rm tot}(s)\dd s\). Its predictable projection onto
\(\mathcal G_t\) has the same density, because that density is already
\(\mathcal G_{t^-}\)-measurable. Equivalently, the tower property applied to
bounded \(\mathcal G_{t^-}\)-predictable test processes gives
\[
 C_{\rm tot}((0,t])-\int_0^t\Lambda_{\rm tot}(s)\dd s
\]
as a \(\mathcal G_t\)-local martingale. Thus, the displayed equation is the
scalar Hawkes intensity equation in the filtration of the unmarked count.
\end{proof}

We next choose the Hawkes mean rate to agree with the collision frequency
from kinetic theory. In dimension three, let \(V=L^3\) and
\(\Sigma=\pi a^2\). The normalized Maxwellian velocity density is
\[
 \phi_T(v)=\left(\frac{m}{2\pi k_{\rm B}T}\right)^{3/2}
 \exp\left(-\frac{m|v|^2}{2k_{\rm B}T}\right)
\]
If \(v,v_*\) are independent
with density \(\phi_T\), then
\begin{equation}\label{eq:relative-speed}
 \bE|v-v_*|=4\sqrt{\frac{k_{\rm B}T}{\pi m}}.
\end{equation}
Indeed, each component of \(v-v_*\) is centered Gaussian with variance
\(2k_{\rm B}T/m\). Thus, \(|v-v_*|\) has the three dimensional Maxwell
distribution with scale \(\sigma=\sqrt{2k_{\rm B}T/m}\), whose mean is
\[
 2\sigma\sqrt{\frac2\pi}
 =4\sqrt{\frac{k_{\rm B}T}{\pi m}}.
\]
Under spatial homogeneity, molecular chaos and the dilute binary collision
approximation, the mean collision rate of one fixed pair is
\begin{equation}\label{eq:kinetic-pair-rate}
 r_{\rm kin}=\frac{\Sigma}{V}\bE|v-v_*|
 =\frac{4\Sigma}{V}\sqrt{\frac{k_{\rm B}T}{\pi m}},
\end{equation}
and the corresponding tagged particle frequency is
\((N-1)r_{\rm kin}\). These are mean rate statements from dilute kinetic
theory
\citep{Cercignani1988,ChapmanCowling1970,CercignaniIllnerPulvirenti1994}.
Indeed, condition on the relative velocity \(g=v-v_*\). During a time
interval \(\dd t\), the
relative center sweeps a collision cylinder of volume
\(\Sigma|g|\dd t\). A homogeneous second center lies in this cylinder with
probability \(\Sigma|g|\dd t/V+o(\dd t)\). Averaging over \(g\) gives the
rate \(\Sigma\bE|g|/V\).

At finite density, a common equilibrium closure uses the
Carnahan--Starling contact correction. The packing fraction is
\[
 \varphi=\frac{\pi Na^3}{6V}
\]
and the associated approximation to the radial distribution at contact is
\begin{equation}\label{eq:carnahan-starling-contact}
 g_{\rm CS}(\varphi)=\frac{1-\varphi/2}{(1-\varphi)^3}.
\end{equation}
We use
\begin{equation}\label{eq:contact-corrected-rate}
 \omega_{\rm CS}=g_{\rm CS}(\varphi)(N-1)r_{\rm kin}
\end{equation}
as the conventional finite density reference value
\citep{CarnahanStarling1969}. Equation
\eqref{eq:contact-corrected-rate} uses an equilibrium closure to approximate
the collision frequency at finite \(N\).

Combining \eqref{eq:mean-rate} and \eqref{eq:kinetic-pair-rate} gives the
calibration
\begin{equation}\label{eq:baseline-calibration}
 \mu=(1-\kappa_G)r_{\rm kin}.
\end{equation}
This choice fixes \(\mu\) by matching the mean rate of a pair. The
higher order statistics remain determined by the chosen Hawkes kernels.

\section{Covariances of the three modes}\label{sec:covariance}

\subsection{Long time covariance}

Assume in this section that the kernels are nonnegative, that
\(\kappa_G<1\), and that
\[
 \int_0^\infty(1+t)h_\alpha(t)\dd t<\infty,\qquad
 \alpha\in\{{\rm same},{\rm share},{\rm dis}\}.
\]
Write \(N(0,t]=(N_e((0,t]))_{e\in\Etwo}\). The standard multivariate Hawkes
covariance formula is
\begin{equation}\label{eq:hawkes-covariance}
\begin{split}
 \Gamma&:=\lim_{t\to\infty}\frac1t\Cov\{N(0,t]\}\\
 &=(I_M-K)^{-1}\diag(r_H\bfone_M)(I_M-K^\top)^{-1}.
\end{split}
\end{equation}
We derive \eqref{eq:hawkes-covariance} by expressing count fluctuations in
terms of the compensated process. Put
\(H(t)=(h_{ef}(t))\), and let
\(R(t)=\sum_{n\geq1}H^{(*n)}(t)\) be the resolvent. The compensated process
\(\dd M(t)=\dd N(t)-\lambda(t)\dd t\) has predictable bracket
\(\diag(\lambda(t))\dd t\), with zero predictable cross brackets because the
coordinates do not jump simultaneously. Put
\[
 G(\dd u)=I_M\delta_0(\dd u)+R(u)\mathbf 1_{\{u>0\}}\dd u.
\]
We first justify the resolvent and the fluctuation representation. Let
\(K_1=\int_0^\infty uH(u)\dd u\), which is finite entrywise by the assumed
kernel moments. For \(n\geq1\), convolution and Fubini's theorem give
\[
 \int_0^\infty H^{(*n)}(u)\dd u=K^n
\]
A second identity is
\[
 \int_0^\infty uH^{(*n)}(u)\dd u
 =\sum_{\ell=0}^{n-1}K^\ell K_1K^{n-1-\ell}.
\]
Choose \(q\) with \(\rho(K)<q<1\). Since the matrices are finite
dimensional, there is a constant \(C_q\) such that
\(\|K^n\|\leq C_qq^n\). The two preceding identities therefore imply
\[
 \sum_{n\geq1}\int_0^\infty(1+u)
       \|H^{(*n)}(u)\|\dd u<\infty,
\]
after replacing the matrix norm by an equivalent entrywise norm if
necessary. Thus \(R\) is a well-defined integrable matrix kernel with a
finite first moment.

The cluster construction gives a bound on second moments that we will use
to justify the limit below. Let \(\mathcal C_f\)
denote the relative marked cluster generated by one immigrant of type \(f\),
and let \(S_f=|\mathcal C_f|\). Subcritical finite type Poisson branching
gives \(\bE S_f^2<\infty\); see
\citep{AthreyaNey1972,Harris1963} and the probability generating function
argument in Remark \ref{rem:scope}. For a deterministic
\(g\in L^2(\bR;\bR^M)\) of compact support, the Poisson integral variance
formula and Cauchy--Schwarz give
\begin{align*}
 &\bE\left|\sum_e\int_{\bR}g_e(s)\dd\widetilde N_e(s)\right|^2\\
 &\quad=\mu\sum_{f\in\Etwo}\int_{\bR}
 \bE\left|\sum_{(\tau,e)\in\mathcal C_f}g_e(u+\tau)\right|^2\dd u\\
 &\quad\leq\mu\sum_{f\in\Etwo}\bE S_f^2\,
 \|g\|_{L^2(\bR;\bR^M)}^2.
\end{align*}
Indeed, for a fixed cluster,
\[
\left|\sum_{(\tau,e)\in\mathcal C_f}g_e(u+\tau)\right|^2
\leq S_f\sum_{(\tau,e)\in\mathcal C_f}|g_e(u+\tau)|^2,
\]
and integration in \(u\) bounds the right side by
\(S_f^2\|g\|_2^2\). Put
\(C_{\rm cl}:=\mu\sum_{f\in\Etwo}\bE S_f^2<\infty\).

Write
\(\dd\widetilde N(t)=\dd N(t)-r_H\bfone_M\dd t\). Subtracting the
mean intensity equation from \eqref{eq:pair-hawkes} gives the stochastic
convolution equation
\[
 \dd\widetilde N=\dd M+H*\dd\widetilde N.
\]
After \(p\) substitutions,
\[
 \dd\widetilde N
 =\left(I_M\delta_0+\sum_{n=1}^{p}H^{(*n)}\right)*\dd M
 +H^{(*(p+1))}*\dd\widetilde N.
\]
To make the passage to the limit explicit, test this identity against
\(g\in L^2(\bR;\bR^M)\) of compact support. The adjoint convolution operator
associated with a matrix kernel \(A\) is
\[
 (T_Ag)(s)=\int_0^\infty A(u)^\top g(s+u)\dd u,
\]
and Young's inequality gives
\[
\|T_Ag\|_2\leq\|A\|_{L^1}\|g\|_2.
\]
The martingale isometry and stationarity therefore show that the
distance in \(L^2\) between two partial martingale convolutions is at most
\[
r_H^{1/2}
\left\|\sum_{n=p+1}^{p'}H^{(*n)}\right\|_{L^1}\|g\|_2,
\]
which tends to zero by the established summability. The preceding cluster
bound gives, for the residual term,
\[
\begin{split}
 &\left\lVert
 \int_{\bR}(T_{H^{(*(p+1))}}g)(s)^\top\dd\widetilde N(s)
 \right\rVert_{L^2}\\
 &\quad\leq C_{\rm cl}^{1/2}\|H^{(*(p+1))}\|_{L^1}\|g\|_2\\
 &\quad\leq C q^{p+1}\|g\|_2.
\end{split}
\]
Thus, the residual tends to zero and the martingale convolutions converge in
\(L^2\). Passing to the limit yields the stationary fluctuation representation,
\[
 \dd N(t)-r_H\bfone_M\dd t
 =\int_{[0,\infty)}G(\dd u)\dd M(t-u).
\]
The measure \(G\) thus describes how a compensated event contributes to
later fluctuations. To count its contribution inside the window \((0,t]\),
define
\[
 \Psi_t(s)=\int_{[0,\infty)}
 \mathbf 1_{\{0<s+u\leq t\}}G(\dd u).
\]
Integration of the preceding representation and the martingale isometry give
\begin{align*}
\Cov\{N(0,t]\}
&=\bE\int_{\bR}\Psi_t(s)\diag\{\lambda(s)\}\Psi_t(s)^\top\dd s\\
&=r_H\int_{\bR}\Psi_t(s)\Psi_t(s)^\top\dd s.
\end{align*}
The second equality uses stationarity and
\(\bE\lambda_e(s)=r_H\). Fubini's theorem yields the exact overlap identity
\[
\begin{split}
 &\frac1t\int_{\bR}\Psi_t(s)\Psi_t(s)^\top\dd s\\
 &\quad=\int_{[0,\infty)^2}
 \left(1-\frac{|u-v|}{t}\right)_+
 G(\dd u)G(\dd v)^\top.
\end{split}
\]
The summability proved above supplies an integrable dominating measure.
Dominated convergence therefore sends the last display to
\[
\left\{\int_{[0,\infty)}G(\dd u)\right\}
\left\{\int_{[0,\infty)}G(\dd u)\right\}^{\!\top}.
\]
Finally,
\[
\int_{[0,\infty)}G(\dd u)
=I_M+\sum_{n\geq1}K^n=(I_M-K)^{-1},
\]
which proves \eqref{eq:hawkes-covariance}. This is also the zero frequency
form of the point process spectral formula
\citep{Bartlett1963,
DaleyVereJones2003,
BacryDelattreHoffmannMuzy2013,
BacryMuzy2016}.

Because \(K\) is symmetric, its three orthogonal modes also separate the
covariance. Applying Theorem \ref{thm:three-mode} to
\eqref{eq:hawkes-covariance} gives
\begin{equation}\label{eq:modal-covariance}
 \Gamma=r_H\left\{
 \frac{P_G}{(1-\kappa_G)^2}
 +\frac{P_P}{(1-\kappa_P)^2}
 +\frac{P_C}{(1-\kappa_C)^2}\right\}.
\end{equation}

\begin{Prop}[Variance factors for observable modes]\label{prop:mode-variance}
Let \(q\neq0\) belong to \(U_a\), where \(a\in\{G,P,C\}\). Then
\begin{equation}\label{eq:q-factor}
 \lim_{t\to\infty}
 \frac{\Var\{q^\top N(0,t]\}}{r_Ht\|q\|^2}
 =\frac1{(1-\kappa_a)^2}.
\end{equation}
For the tagged count \(C_i(0,t]\),
\begin{equation}\label{eq:tagged-factor}
 \lim_{t\to\infty}\frac{\Var\{C_i(0,t]\}}
 {(N-1)r_Ht}
 =\frac2N\frac1{(1-\kappa_G)^2}
 +\frac{N-2}{N}\frac1{(1-\kappa_P)^2}.
\end{equation}
Every contrast \(C_i-C_j\) belongs to \(U_P\), while every \(q\in U_C\)
satisfies \(Bq=0\) and is therefore absent from all tagged counts.
\end{Prop}

\begin{proof}
If \(q\in U_a\), then \(P_aq=q\) and the other two projectors annihilate
\(q\). Equation \eqref{eq:q-factor} follows by applying \(q^\top(\cdot)q\)
to \eqref{eq:modal-covariance}.

Let \(b_i=B^\top e_i\), so that \(C_i=b_i^\top N\). Its global projection is
\[
 P_Gb_i=\frac{\bfone_M^\top b_i}{M}\bfone_M
 =\frac{N-1}{M}\bfone_M=\frac2N\bfone_M,
\]
because \(b_i\) is the indicator of the \(N-1\) edges incident to particle
\(i\). Moreover, \(b_i\in\operatorname{range}(B^\top)\), and hence
\(P_Cb_i=0\) and \(P_Pb_i=b_i-P_Gb_i\). Orthogonality gives
\[
 \|P_Gb_i\|^2
 =\frac4{N^2}M=\frac{2(N-1)}N
\]
and, since \(\|b_i\|^2=N-1\),
\[
 \|P_Pb_i\|^2
 =\|b_i\|^2-\|P_Gb_i\|^2
 =\frac{(N-1)(N-2)}N.
\]
It follows from \eqref{eq:modal-covariance} that
\[
\begin{split}
 &\lim_{t\to\infty}\frac1{r_Ht}\Var\{C_i(0,t]\}\\
 &\quad=\frac{2(N-1)/N}{(1-\kappa_G)^2}
 +\frac{(N-1)(N-2)/N}{(1-\kappa_P)^2}.
\end{split}
\]
Division by \(N-1\) gives \eqref{eq:tagged-factor}. Finally,
\(b_i-b_j=B^\top(e_i-e_j)\in U_P\), and \(Bq=0\) is the definition of
\(U_C\).
\end{proof}

\subsection{Exact finite-window factors for exponential kernels}

For exponential kernels, the covariance can also be calculated at each
finite observation window. Suppose that the three kernels have a common
exponential
basis:
\begin{equation}\label{eq:exponential-kernels}
 h_\alpha(t)=\beta k_\alpha e^{-\beta t}\mathbf 1_{\{t\geq0\}},
 \qquad \beta>0.
\end{equation}
For \(q\in U_a\), define the normalized variance factor of width \(\Delta\),
\[
 F_a(\Delta)=
 \frac{\Var\{q^\top N(0,\Delta]\}}
 {r_H\Delta\|q\|^2}.
\]
The modal spectral density divided by \(r_H\) is
\[
 \left|1-\frac{\kappa_a\beta}{\beta+i\omega}\right|^{-2}
 =\frac{\beta^2+\omega^2}
 {\beta^2(1-\kappa_a)^2+\omega^2}.
\]
Because \(K\) is symmetric and \(\rho(K)=\kappa_G<1\), every modal
eigenvalue satisfies \(|\kappa_a|<1\). Hence
\(a_a=\beta(1-\kappa_a)>0\). Then
\begin{align*}
\frac{\beta^2+\omega^2}{a_a^2+\omega^2}
&=1+\frac{\beta^2-a_a^2}{a_a^2+\omega^2}\\
&=1+\frac{\beta^2\kappa_a(2-\kappa_a)}
{a_a^2+\omega^2}.
\end{align*}
Under the convention
\(\cF^{-1}f(u)=(2\pi)^{-1}\int_{\bR}e^{i\omega u}f(\omega)\dd\omega\),
\[
 \cF^{-1}\left\{\frac1{a_a^2+\omega^2}\right\}(u)
 =\frac{e^{-a_a|u|}}{2a_a}.
\]
The inverse Fourier transform therefore gives the normalized covariance
measure,
\[
 \delta_0(\dd u)+
 c_a e^{-a_a|u|}\dd u,\qquad
 c_a=\frac{\beta\kappa_a(2-\kappa_a)}{2(1-\kappa_a)}.
\]
Integrating this measure over pairs of times in \([0,\Delta]\) gives
\begin{align*}
F_a(\Delta)
&=1+\frac{c_a}{\Delta}
\int_0^\Delta\int_0^\Delta e^{-a_a|u-v|}\dd u\dd v\\
&=1+\frac{2c_a}{\Delta}
\int_0^\Delta(\Delta-u)e^{-a_au}\dd u.
\end{align*}
The integral in the last expression is
\[
 \int_0^\Delta(\Delta-u)e^{-a_au}\dd u
 =\frac{\Delta}{a_a}-\frac{1-e^{-a_a\Delta}}{a_a^2}.
\]
Consequently
\begin{align*}
F_a(\Delta)
&=1+\frac{2c_a}{a_a}
 -\frac{2c_a\{1-e^{-a_a\Delta}\}}{a_a^2\Delta}\\
&=\frac1{(1-\kappa_a)^2}
 -\frac{\kappa_a(2-\kappa_a)
 \{1-e^{-\beta(1-\kappa_a)\Delta}\}}
 {\beta(1-\kappa_a)^3\Delta},
\end{align*}
where we have used the identity
\[
1+\frac{2c_a}{a_a}
=1+\frac{\kappa_a(2-\kappa_a)}{(1-\kappa_a)^2}
=\frac1{(1-\kappa_a)^2}.
\]
Thus,
\begin{equation}\label{eq:finite-window-factor}
 F_a(\Delta)=\frac1{(1-\kappa_a)^2}
 -\frac{\kappa_a(2-\kappa_a)
 \{1-e^{-\beta(1-\kappa_a)\Delta}\}}
 {\beta(1-\kappa_a)^3\Delta}.
\end{equation}
In particular, \(F_a(\Delta)\to(1-\kappa_a)^{-2}\) as
\(\Delta\to\infty\). The second term in
Eq.~\eqref{eq:finite-window-factor} gives the difference between a finite
window and this limit. Section~\ref{sec:numerical} compares the formula
with numerical estimates.

\section{Numerical results}\label{sec:numerical}

We first compare the covariance formulas with fluctuations measured in
pair-Hawkes simulations. We then examine the collision counts of hard
particles in two and three dimensions, using their modal variances and
cumulants to describe the fluctuations. The numerical data and
reproducibility materials are available in
\cite{HernandezRuiz2026Data}. The results below concern finite systems and
finite observation windows.

\subsection{Fluctuations of the pair-Hawkes process}

For exponential kernels, the three gains determine how strongly each mode
fluctuates. We use
\[
 \begin{gathered}
 (\kappa_G,\kappa_P,\kappa_C)=(0.45,0.22,0.06),\\
 \beta=2,\qquad r_H=0.006.
 \end{gathered}
\]
For each particle number, Eqs.~\eqref{eq:inverse-same}--\eqref{eq:inverse-dis}
give the kernel masses that preserve these gains. Let
\(\widehat\Gamma_\Delta\) denote the empirical covariance of the pair
counts over a window of width \(\Delta\). Its normalized projection is
\begin{equation}\label{eq:mode-estimator}
 \widehat F_a(\Delta)
 =\frac{\tr(P_a\widehat\Gamma_\Delta)}{d_a r_H\Delta},
 \qquad d_a=\dim U_a,
\end{equation}
for \(a\in\{G,P,C\}\). The trace sums the variances over an orthonormal
basis of the mode. Division by \(d_a\) gives the average variance per
direction.
Division by \(r_H\Delta\) then expresses it relative to the Poisson scale.

\begin{figure*}[t]
\centering
\includegraphics[width=\textwidth]{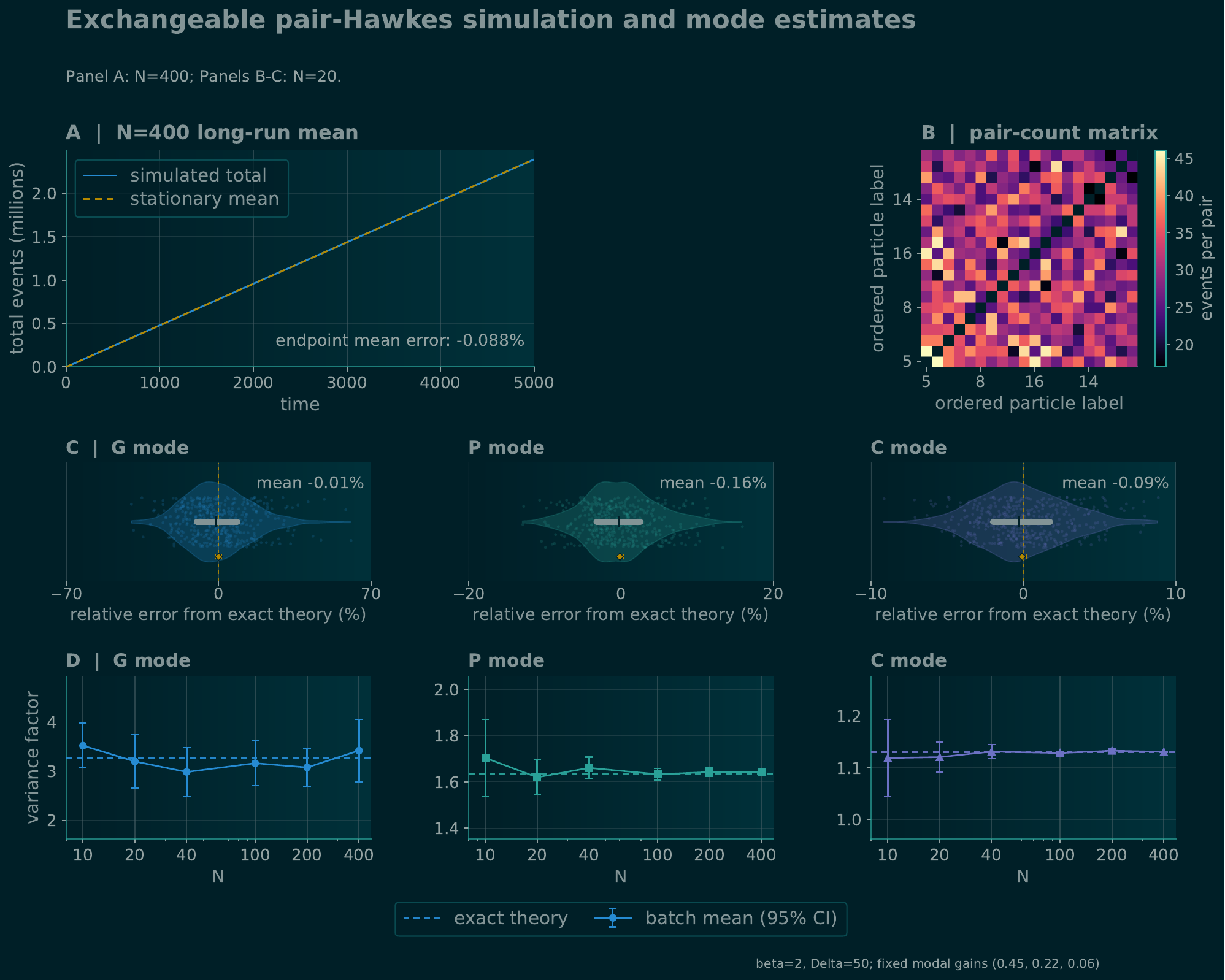}
\caption{Fluctuations of the exchangeable pair-Hawkes process. Panel A
compares the cumulative total with its stationary mean. Panel B shows pair
counts with particle labels ordered by their tagged counts. Panel C
displays relative errors against Eq.~\eqref{eq:finite-window-factor}:
dots and violins show the empirical distributions, thick segments the
interquartile ranges, ticks the medians, and diamonds with bars the means
and their 95\% confidence intervals. The three strips have separate
horizontal scales. Panel D compares the mode factors across particle
numbers. Error bars are pointwise 95\% intervals across batch estimates;
dashed lines give the exact finite-window factors.}
\label{fig:simulation}
\end{figure*}

Figure~\ref{fig:simulation} shows the cumulative count, the pair counts
and the three modal estimates. At \(N=20\) and \(\Delta=50\), the mean
estimates in Table~\ref{tab:mode-validation} differ from the exact
finite-window factors by at most \(0.16\%\). The global mode has the
largest factor and the greatest dispersion among the estimates.

\begin{table*}[t]
\caption{Pair-Hawkes mode factors at \(N=20\) and \(\Delta=50\).
Intervals describe the uncertainty in the mean across independent
realizations. The exact values are given by
Eq.~\eqref{eq:finite-window-factor}.}
\label{tab:mode-validation}
\centering
\setlength{\tabcolsep}{6pt}
\begin{tabular}{lccc}
\toprule
Mode & Mean estimate & 95\% confidence interval & Exact factor\\
\midrule
Global & 3.2636 & [3.2148, 3.3125] & 3.2639\\
Tagged particle & 1.6328 & [1.6254, 1.6402] & 1.6354\\
Cycle & 1.1293 & [1.1261, 1.1325] & 1.1303\\
\bottomrule
\end{tabular}
\end{table*}

A separate comparison varies the particle number while keeping the gains
fixed. The same theoretical factors therefore apply at every \(N\).
Table~\ref{tab:finite-size-validation} gives the corresponding estimates.
Every pointwise interval in
Fig.~\ref{fig:simulation}(D) contains its corresponding exact factor.
As the tagged particle and cycle spaces grow in dimension, their estimates
become less dispersed. The global space remains one dimensional. These
observations support the covariance formulas for this parametrization of
the reference process. The finite-size behavior of a physical gas requires
a separate analysis.

\begin{table}[t]
\caption{Mean mode-factor estimates at different particle numbers, with
\(\Delta=50\) and fixed modal gains. Figure~\ref{fig:simulation}(D)
shows the corresponding pointwise 95\% intervals.}
\label{tab:finite-size-validation}
\centering
\setlength{\tabcolsep}{6pt}
\begin{tabular}{rccc}
\toprule
\(N\) & \(\widehat F_G\) & \(\widehat F_P\) & \(\widehat F_C\)\\
\midrule
10  & 3.5223 & 1.7032 & 1.1190\\
20  & 3.2003 & 1.6195 & 1.1207\\
40  & 2.9845 & 1.6602 & 1.1313\\
100 & 3.1638 & 1.6325 & 1.1286\\
200 & 3.0772 & 1.6420 & 1.1330\\
400 & 3.4215 & 1.6406 & 1.1308\\
\midrule
Exact & 3.2639 & 1.6354 & 1.1303\\
\bottomrule
\end{tabular}
\end{table}

The distance between the cumulative count and its mean in
Fig.~\ref{fig:simulation}(A) can be read on the scale of the global
fluctuations. For \(S(t)=\sum_{e\in\Etwo}N_e(0,t]\), the identity
\(P_G\bfone_M=\bfone_M\) gives
\begin{equation}\label{eq:total-count-finite-variance}
 \Var S(t)=M r_H t F_G(t).
\end{equation}
For large \(t\), the absolute fluctuation scale grows as \(t^{1/2}\),
while the mean
grows linearly in time. The relative fluctuation therefore decays as
\(t^{-1/2}\). At the endpoint of the plotted record, the difference from
the stationary mean is \(-0.75\) standard deviations according to this
formula.

\begin{figure*}[t]
\centering
\includegraphics[width=\textwidth]{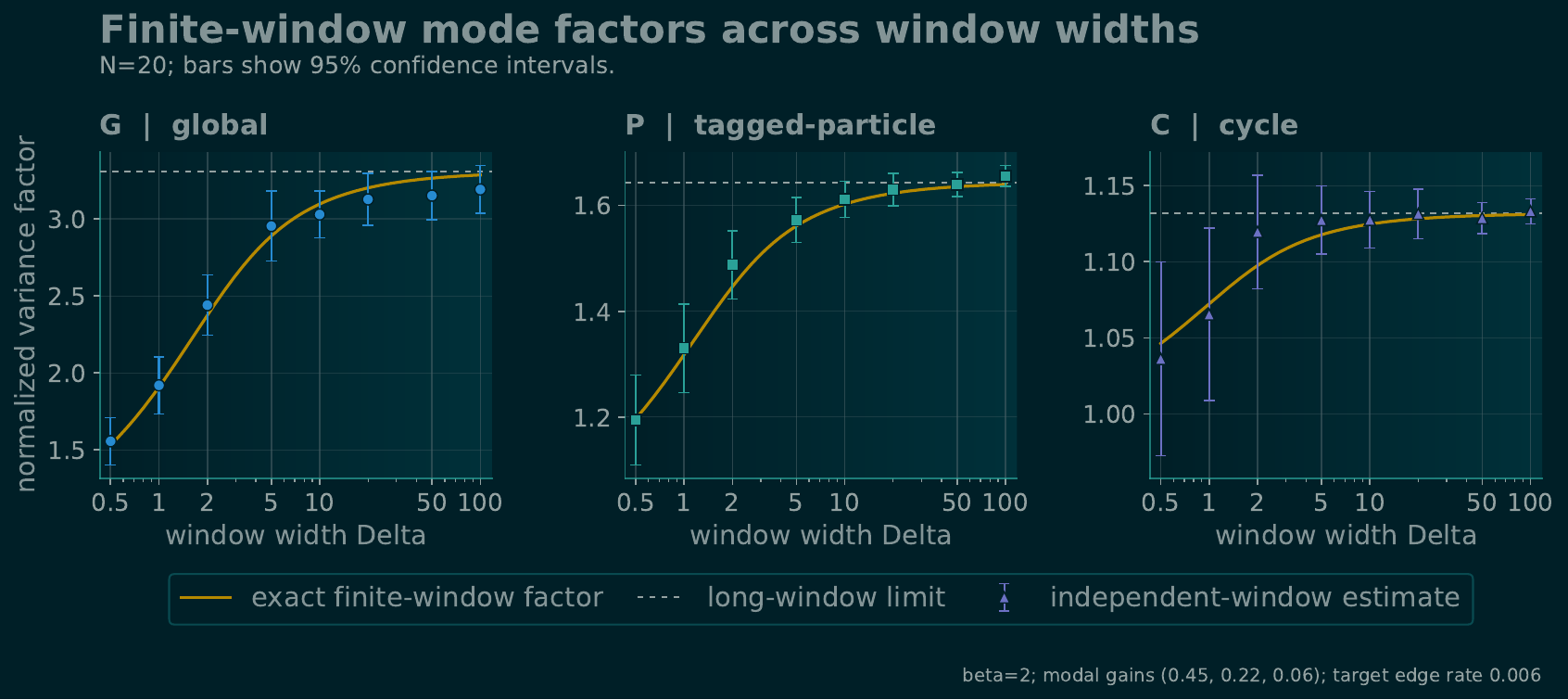}
\caption{Dependence of the mode factors on the observation-window width
\(\Delta\), for \(N=20\). Points are mean estimates with pointwise
95\% confidence intervals across batch estimates. Solid curves give
Eq.~\eqref{eq:finite-window-factor}, and horizontal dashed lines give
the long-time limits \((1-\kappa_a)^{-2}\).}
\label{fig:window-sweep}
\end{figure*}

Figure~\ref{fig:window-sweep} follows the factors as the observation window
widens. Each estimate is compared with the full expression in
Eq.~\eqref{eq:finite-window-factor}. Every displayed pointwise 95\%
interval contains the corresponding exact value, and the largest absolute
standardized discrepancy is \(1.577\). As \(\Delta\) increases, the
factors approach \((1-\kappa_a)^{-2}\). Their separation from these limits
at shorter windows shows how much the finite-window term contributes.

\subsection{Tagged collision counts in two dimensions}
\label{subsec:hard-disk-check}

For the hard-particle comparisons, write \(\widehat c_j(\Delta)\) for
the tagged-count cumulants per unit time, estimated over a finite
window. We reserve
\(c_j\) for long-time rates and \(c_j^{\rm S}\) for the first Sonine
rates. To compare observation windows on the collision time scale, we use
\(\theta=\widehat\omega\Delta\), where \(\widehat\omega\) is the measured
tagged particle rate. This expresses the window width in units of the mean
time between tagged collisions.

\begin{figure*}[t]
\centering
\includegraphics[width=\textwidth]{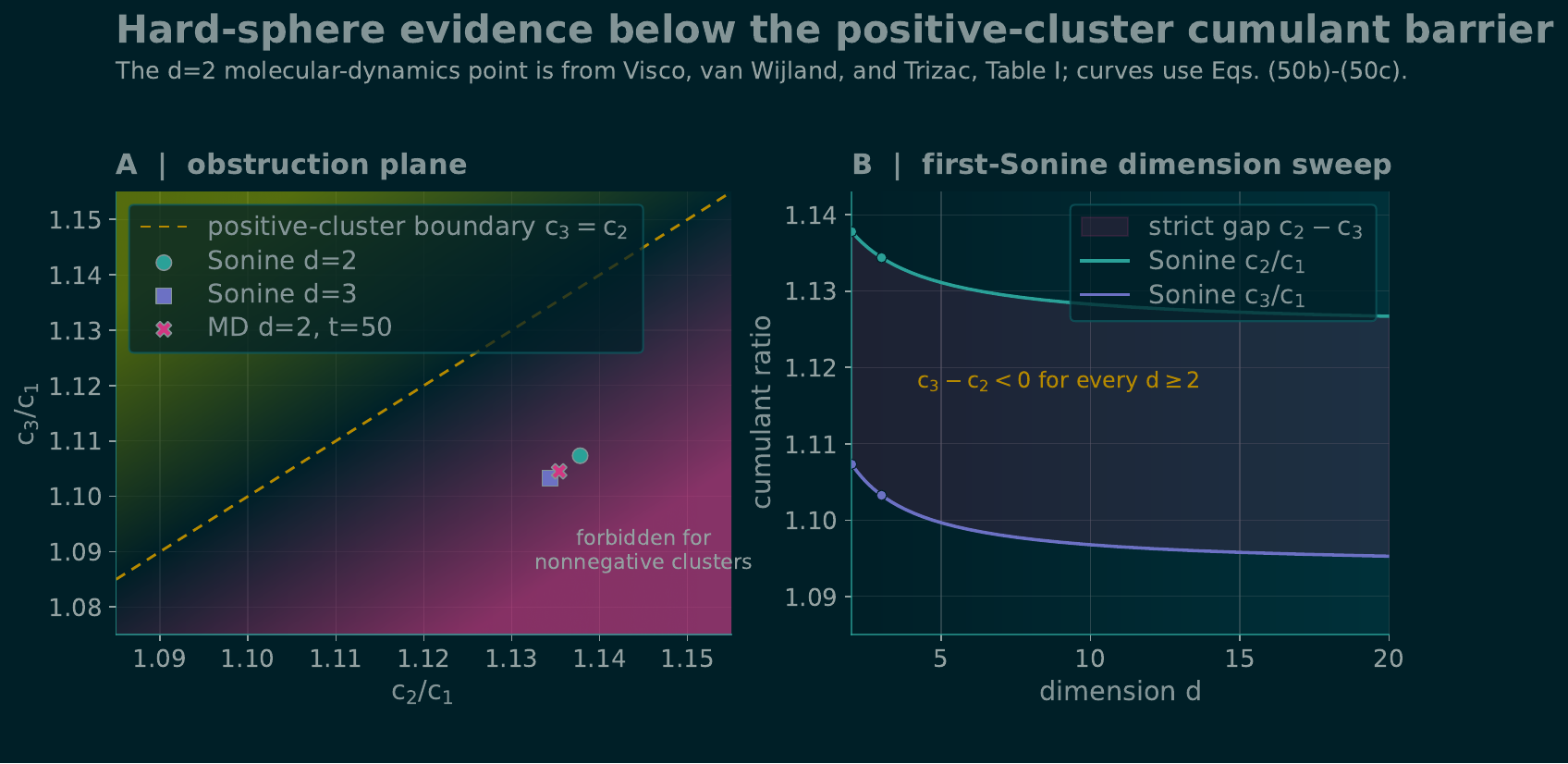}
\caption{Hard-particle cumulants and the positive cluster inequality.
Panel A shows the first Sonine values in dimensions two and three and the
published two-dimensional molecular dynamics estimate at \(\omega t=50\)
\cite{ViscoVanWijlandTrizac2008}. All lie below the diagonal corresponding
to equal second and third cumulants. Panel B compares the first Sonine
ratios as the dimension varies. Their separation gives the strict gap in
Eq.~\eqref{eq:sonine-gap}. The molecular dynamics point describes a
finite observation window; the Sonine values are long-time coefficients
within the kinetic approximation.}
\label{fig:cumulant-evidence}
\end{figure*}

Figure~\ref{fig:cumulant-evidence} compares the first Sonine ratios with
the molecular dynamics point reported by
\citet{ViscoVanWijlandTrizac2008} at \(\omega t=50\). The first Sonine
ratios \(c_2^{\rm S}/c_1^{\rm S}\) and
\(c_3^{\rm S}/c_1^{\rm S}\) are, respectively,
\[
 \begin{array}{c|cc}
 d&c_2^{\rm S}/c_1^{\rm S}&c_3^{\rm S}/c_1^{\rm S}\\ \hline
 2&1.1377841&1.1073257\\
 3&1.1343750&1.1032510.
 \end{array}
\]
The published planar molecular dynamics ratios are \(1.1354\) and
\(1.1045\) at \(\omega t=50\), giving a normalized difference of
\(-0.0309\). At \(\omega t=10\), the reported ratios are \(1.1228\)
and \(1.1282\), whose difference is \(0.0054\). The difference therefore
changes sign between the two reported windows. The first Sonine gap in
Eq.~\eqref{eq:sonine-gap} describes the long-time kinetic approximation.
The molecular dynamics points show the ordering at the reported finite
windows.

\begin{figure*}[t]
\centering
\includegraphics[width=\textwidth]{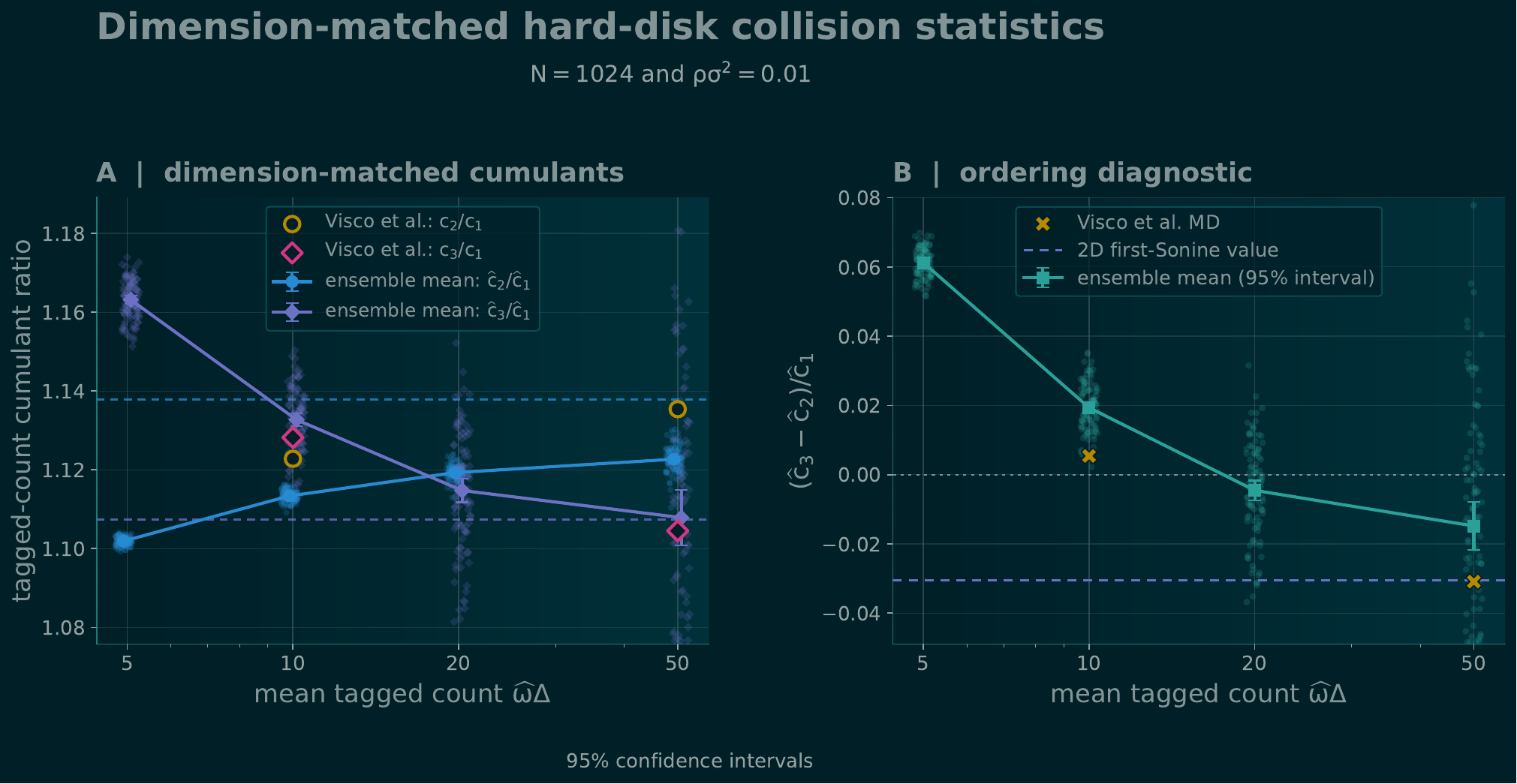}
\caption{Tagged collision-count fluctuations in two dimensions. Panel A
compares the mean cumulant-ratio estimates with the molecular dynamics
results of \citet{ViscoVanWijlandTrizac2008}; dashed lines give the
two-dimensional first Sonine values. Panel B shows the normalized
difference between the third and second cumulants. Faint points represent
individual realizations, and error bars are pointwise 95\% Studentized
intervals across realization estimates. At the longest window, the
interval for the mean difference lies below zero. The published points
have no error bars because the source tabulates no uncertainty estimates.}
\label{fig:hard-disk-2d}
\end{figure*}

The measured planar tagged particle rate is
\(\widehat\omega=0.0439391\), with a 95\% interval for the mean across
independent records of \([0.0439347,0.0439437]\). The intervals below are
descriptive summaries for means across the independent planar records.

Our independent planar calculation gives the ratios in
Table~\ref{tab:planar-cumulants} and Fig.~\ref{fig:hard-disk-2d}.
The mean normalized difference changes from positive at \(\theta=10\) to
negative at \(\theta=50\). At the longer window, the entire interval
\([-0.0217,-0.0079]\) lies below zero. The planar records therefore support
a mean third cumulant ratio smaller than the second at that finite window.
The published values show the same change in ordering. The numerical
ratios differ: at \(\theta=50\), our second and third
ratios differ from the published values by \(-0.0127\) and \(0.0034\),
respectively. The asymptotic sign for exact hard-particle dynamics remains
open here.

\begin{table*}[t]
\caption{Finite-window tagged collision cumulants in two dimensions.
Brackets give pointwise 95\% Studentized intervals for means across
independent realizations. The difference changes sign between the two
windows.}
\label{tab:planar-cumulants}
\centering
\setlength{\tabcolsep}{6pt}
\begin{tabular}{cccc}
\toprule
\(\theta\) & \(\widehat c_2/\widehat c_1\)
& \(\widehat c_3/\widehat c_1\)
& \((\widehat c_3-\widehat c_2)/\widehat c_1\)\\
\midrule
10 & 1.1134 [1.1131, 1.1137] & 1.1327 [1.1311, 1.1343]
& 0.0193 [0.0179, 0.0208]\\
50 & 1.1227 [1.1220, 1.1234] & 1.1079 [1.1008, 1.1149]
& $-0.0148$ [$-0.0217$, $-0.0079$]\\
\bottomrule
\end{tabular}
\end{table*}

\subsection{Collision statistics in three dimensions}
\label{subsec:event-driven-check}

\begin{figure*}[t]
\centering
\includegraphics[width=\textwidth]{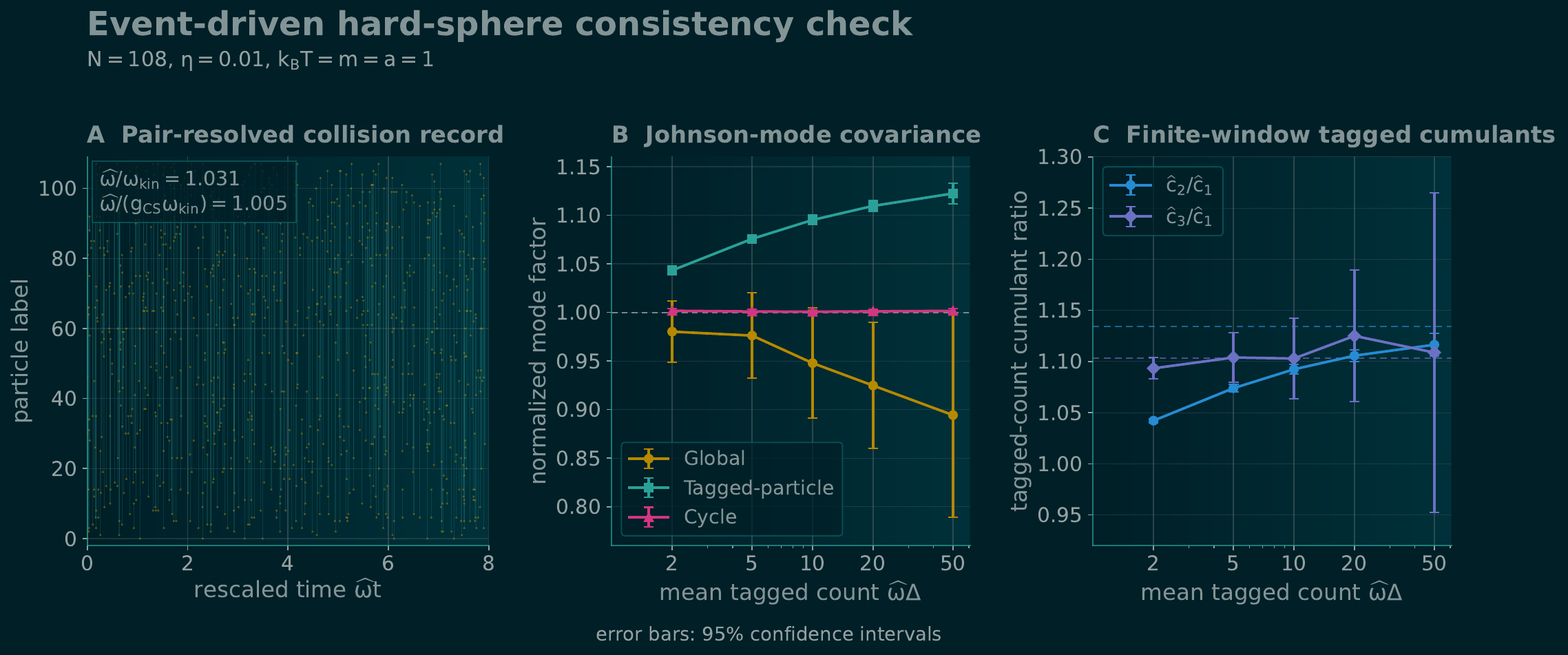}
\caption{Collision statistics of a finite three-dimensional hard-sphere
system. Panel A shows the pair-resolved record and compares the measured
collision rate with the dilute and contact-corrected kinetic values.
Panel B gives the normalized modal projections in
Eq.~\eqref{eq:hard-sphere-mode-factor}. Panel C shows the tagged cumulant
ratios, with dashed lines marking the three-dimensional first Sonine
values. Error bars in B and C are pointwise 95\% Studentized intervals
from within-record batch means. The cumulant uncertainty leaves their
ordering unresolved.}
\label{fig:hard-sphere-validation}
\end{figure*}

The three-dimensional calculation gives a measured tagged particle rate
\[
 \widehat\omega=0.1383176\quad[0.1380166,0.1386187].
\]
The bracket is a descriptive 95\% interval from within-record batch
means. The dilute and contact-corrected kinetic rates from
Eqs.~\eqref{eq:kinetic-pair-rate} and
\eqref{eq:contact-corrected-rate} are
\[
 \omega_{\rm kin}=0.1341517,\qquad
 \omega_{\rm CS}=0.1375669.
\]
The measured rate is \(3.11\%\) above the dilute value and \(0.55\%\)
above the contact-corrected value. Figure~\ref{fig:hard-sphere-validation}(A)
shows the pair-resolved collision record alongside these rates. The
measured rate is closer to the contact-corrected value for this finite
system.

We can also use the Johnson projectors to separate the fluctuations of
the physical pair counts into the same three modes. For their empirical
covariance \(\widehat\Gamma_\Delta^{\rm HS}\) and mean rate per pair
\(\widehat r_e\), define
\begin{equation}\label{eq:hard-sphere-mode-factor}
 \widehat F_a^{\rm HS}(\Delta)
 =\frac{\tr(P_a\widehat\Gamma_\Delta^{\rm HS})}
 {d_a\widehat r_e\Delta},\qquad a\in\{G,P,C\}.
\end{equation}
These are normalized projections of the measured covariance.
 At \(\theta=50\), they give the values in
Table~\ref{tab:physical-modes}. Relative to the Poisson scale, the global
factor lies below one, the tagged particle factor above one, and the cycle
factor is close to one. Figure~\ref{fig:hard-sphere-validation}(B) shows how
these factors vary with the window width in the finite physical system.

\begin{table*}[t]
\caption{Normalized covariance projections for the three-dimensional
collision record at \(\theta=50\). Brackets give descriptive pointwise
95\% Studentized intervals from within-record batch means.}
\label{tab:physical-modes}
\centering
\setlength{\tabcolsep}{6pt}
\begin{tabular}{lcc}
\toprule
Mode & \(\widehat F_a^{\rm HS}\) & 95\% interval\\
\midrule
Global & 0.8943 & [0.7894, 0.9992]\\
Tagged particle & 1.1224 & [1.1118, 1.1330]\\
Cycle & 1.0015 & [0.9990, 1.0041]\\
\bottomrule
\end{tabular}
\end{table*}

At the same window, the tagged cumulant ratios in
Fig.~\ref{fig:hard-sphere-validation}(C) are
\begin{align*}
 \widehat c_2/\widehat c_1&=1.1164\quad[1.1053,1.1274],\\
 \widehat c_3/\widehat c_1&=1.1087\quad[0.9527,1.2648].
\end{align*}
Their normalized difference is
\[
 \frac{\widehat c_3-\widehat c_2}{\widehat c_1}
 =-0.0076\quad[-0.1612,0.1459].
\]
This interval contains both positive and negative values, leaving the
ordering unresolved for the three-dimensional record. Its intervals
summarize variation within the trajectory. The independent planar
realizations provide the finite-window comparison described above. The
analytic obstruction in Sec.~\ref{sec:obstruction} follows from the first
Sonine rates.

\clearpage

\section{Latent state and signed extensions}\label{sec:extensions}

The proof of the cumulant inequality uses independent Poisson immigrant
clusters with nonnegative integer contributions. We consider two extensions
with different sources of dependence between
events. In a latent state model, the
state left by an event determines the next waiting time. Hard-sphere
scattering leads to the linear Boltzmann operator within the kinetic
description, and a separate two-state example gives \(c_3<c_2\). In a signed
Hawkes model, a past event may reduce the present intensity. We prove a
sufficient stationarity condition for its clipped form. Its cumulants remain
to be determined.

\subsection{A normalized latent kinetic state}

Let \(W(v'|v)\) be a nonnegative measurable transition rate density.
Conditional on the present velocity \(v\), the probability of a collision
that moves the tagged velocity into \(\dd v'\) during \(\dd t\) is
\[
 W(v'|v)\dd v'\dd t+o(\dd t).
\]
The total collision hazard and the conditional postcollision density are
\begin{equation}\label{eq:kinetic-rate-transition}
 r(v)=\int_{\bR^d}W(v'|v)\dd v',\qquad
 p(v'|v)=\frac{W(v'|v)}{r(v)}
\end{equation}
whenever \(r(v)>0\). Let \(\phi(v)\) be a normalized equilibrium velocity
density, \(\int_{\bR^d}\phi(v)\dd v=1\). Assume that the associated pure jump
process is nonexplosive, that \(r(v)<\infty\) for \(\phi\) almost every \(v\),
and that \(\phi\) is invariant.
\[
 \int_{\bR^d}\phi(v)W(v'|v)\dd v=r(v')\phi(v')
 \quad\text{for almost every \(v'\)}.
\]
Assume also that the mean stationary collision frequency satisfies
\begin{equation}\label{eq:omega-normalization}
 0<\omega=\int_{\bR^d}r(v)\phi(v)\dd v<\infty.
\end{equation}
At stationarity, the mean number of jumps per unit time whose velocity before
a jump lies in \(\dd v\) is \(r(v)\phi(v)\dd v\). Normalizing this Campbell
measure shows that sampling the precollision velocity at an event time gives
the Palm density
\begin{equation}\label{eq:collision-biased-density}
 \phi_{\rm coll}(v)=\frac{r(v)\phi(v)}{\omega}.
\end{equation}
Sampling at collision times weights the equilibrium velocity density by
the hazard. The factor \(\omega\) in \eqref{eq:collision-biased-density}
normalizes this collision-biased distribution to have total mass one.

When every jump is counted, exponential tilting multiplies each jump
contribution by \(e^s\), giving the operator
\begin{equation}\label{eq:tilted-generator}
 (\cL_s f)(v)=\int_{\bR^d}W(v'|v)
 \{e^s f(v')-f(v)\}\dd v'.
\end{equation}
Under the following spectral assumptions, the dominant eigenvalue
determines the long-time cumulants. Let
\(\mathcal X\) be a Banach lattice, with complexification
\(\mathcal X_{\mathbb C}\), containing \(\bfone\). Suppose that, for complex
\(|s|<\varepsilon\), the operators \(\cL_s\) have a common dense domain and
form an analytic family. Assume that the tilted strongly continuous
semigroups exist throughout this complex disc and depend analytically on
\(s\). Assume directly that there are analytic functions
\(\alpha(s)\), \(r_s\) and \(\ell_s\), and constants \(C,\gamma>0\), such that,
locally uniformly,
\[
 e^{t\cL_s}
 =e^{t\alpha(s)}r_s\otimes\ell_s+\mathcal R_s(t),\qquad
 \|\mathcal R_s(t)\|
 \leq Ce^{t\{\Re\alpha(s)-\gamma\}},
\]
where \(\langle\ell_s,r_s\rangle=1\). Let the initial probability law
\(\nu\) define a continuous functional on \(\mathcal X\), normalized by
\(\langle\nu,\bfone\rangle=1\), and assume that
\[
 a(s):=\langle\nu,r_s\rangle\langle\ell_s,\bfone\rangle
\]
is bounded away from zero on a possibly smaller complex disc. At \(s=0\),
we take \(\alpha(0)=0\), \(r_0=\bfone\), and
\(\langle\ell_0,\bfone\rangle=1\), so that \(a(0)=1\). These assumptions
justify locally uniform convergence of the scaled
cumulant generating functions and convergence of their derivatives. The next
proposition gives two settings
where conditions on the jump process imply the uniform remainder bound and
the nonvanishing coefficient.

\begin{Prop}[Sufficient conditions for the kinetic spectral hypothesis]
\label{prop:kinetic-spectral}
Interpret \eqref{eq:tilted-generator} on a measurable velocity space \(E\),
with integrals replaced by sums when \(E\) is finite. The spectral assumptions
above hold for every initial probability law \(\nu\) in either of the
following cases:
\begin{enumerate}
\item \(E\) is finite and the conservative generator \(\cL_0\) is
irreducible;
\item \(E=\bR^d\), and there are constants \(0<r_-\leq r_+<\infty\) for which
the following bound holds.
\[
 r_-\leq r(v)\leq r_+,
\]
and, for some integer \(m\geq1\), some \(\epsilon>0\), and some probability
measure \(\eta\), the embedded transition kernel satisfies the following
minorization after \(m\) steps,
\[
 p^m(A|v)\geq\epsilon\eta(A)
\]
for every Borel set \(A\) and every \(v\).
\end{enumerate}
One may take \(\mathcal X=\bR^E\) in the first case and
\(\mathcal X=\mathcal B_b(E;\bR)\), with the supremum norm, in the second.
In both cases the stationary law is unique. In the second case it has a
density whenever \(W(v'|v)\) is a transition rate density.
\end{Prop}

\begin{proof}
In the finite-state case, irreducibility of the conservative generator gives
a unique invariant probability \(\pi\), with
\(\ker\cL_0=\operatorname{span}\{\bfone\}\). Every other eigenvalue has
strictly negative real part. Hence, on
\mbox{\(\mathcal X_0=\{f:\pi(f)=0\}\)}, the finite matrix
exponential satisfies
\(\|e^{t\cL_0}|_{\mathcal X_0}\|\leq C_0e^{-\gamma_0t}\) for suitable
\(C_0,\gamma_0>0\); any polynomial factors from Jordan blocks are absorbed
by a smaller exponential rate. Thus, zero is a simple isolated eigenvalue. The
associated decomposition is
\[
\begin{gathered}
 e^{t\cL_0}=\Pi_0+\mathcal R_0(t),\qquad
 \Pi_0f=\pi(f)\bfone,\\
 \|\mathcal R_0(t)\|\leq C_0e^{-\gamma_0t}.
\end{gathered}
\]
The finitely many exit rates are bounded, so the chain is nonexplosive.
The finite-state tilted generators are matrices whose entries are entire
functions of \(s\).

Consider the second case. The rate bounds imply
\(\|\cL_0f\|_\infty\leq2r_+\|f\|_\infty\), so \(\cL_0\) is bounded and
generates a uniformly continuous Markov semigroup \(P_t\) on
\(\mathcal B_b(E)\). For ordered jump times
\(0<t_1<\cdots<t_m<t\), the product of the \(m\) jump rates is at least
\(r_-^m\), while the product of all holding time survival factors is at
least \(e^{-r_+t}\). Integrating the jump times over their simplex and the
successive states through \(p\) gives
\[
\begin{split}
 P_t(v,A)
 &\geq e^{-r_+t}r_-^m
 \int_{0<t_1<\cdots<t_m<t}\dd t_1\cdots\dd t_m\,
 p^m(A|v)\\
 &=e^{-r_+t}\frac{(r_-t)^m}{m!}\,p^m(A|v)
 \geq\delta_t\eta(A),
\end{split}
\]
where
\[
 \delta_t:=\epsilon e^{-r_+t}\frac{(r_-t)^m}{m!}>0.
\]

Fix \(t_0>0\), put \(P=P_{t_0}\) and \(\delta=\delta_{t_0}\). The
minorization permits the decomposition
\[
 P(v,\cdot)=\delta\eta(\cdot)+(1-\delta)Q(v,\cdot)
\]
for a Markov kernel \(Q\). Therefore, for probability measures
\(\mu,\mu'\),
\[
 \|\mu P-\mu'P\|_{\rm TV}
 \leq(1-\delta)\|\mu-\mu'\|_{\rm TV}.
\]
In particular,
\[
 \|\mu P^{n+1}-\mu P^n\|_{\rm TV}
 \leq(1-\delta)^n\|\mu P-\mu\|_{\rm TV},
\]
so \((\mu P^n)\) is Cauchy in total variation. Its limit \(\pi\) is invariant,
because the map \(\mu\mapsto\mu P\) is continuous in total variation. The
same contraction proves uniqueness and gives
\[
 \sup_{v\in E}\|P^n(v,\cdot)-\pi\|_{\rm TV}
 \leq2(1-\delta)^n.
\]
For \(u\geq0\), the measure \(\pi P_u\) is also invariant for \(P=P_{t_0}\),
because \(P_u\) and \(P_{t_0}\) commute. Uniqueness therefore gives
\(\pi P_u=\pi\). The semigroup property, applied to
\(t=nt_0+u\), \(0\leq u<t_0\), then gives constants
\(C_0,\gamma_0>0\) such that
\[
 \|P_t-\Pi_0\|_{\mathcal B_b(E)\to\mathcal B_b(E)}
 \leq C_0e^{-\gamma_0t},\qquad
 \Pi_0f=\pi(f)\bfone.
\]
This is the whole space Doeblin argument; see also
\citep{MeynTweedie2009}.

The last bound gives the spectral gap needed for the generator. The
invariant closed subspace
\mbox{\(\mathcal X_0=\{f:\pi(f)=0\}\)} satisfies
\(\|P_t|_{\mathcal X_0}\|\leq C_0e^{-\gamma_0t}\). Hence, for
\(\Re z>-\gamma_0\), the Laplace integral,
\[
 (z-\cL_0|_{\mathcal X_0})^{-1}f
 =\int_0^\infty e^{-zt}P_tf\dd t
\]
converges in operator norm. Thus, the spectrum of
\(\cL_0|_{\mathcal X_0}\) lies in
\(\{\Re z\leq-\gamma_0\}\), while \(\cL_0\bfone=0\). It follows that zero is
a simple isolated spectral point of \(\cL_0\).

If \(W\) is a transition rate density, invariance
\(\int_E\cL_0f(v)\pi(\dd v)=0\), first for bounded indicator functions and
then by a monotone class argument, identifies the measures,
\[
 r(v)\pi(\dd v)
 =\left\{\int_E W(v|u)\pi(\dd u)\right\}\dd v.
\]
Since \(r(v)\geq r_->0\), division by \(r(v)\) shows that \(\pi\) has a
density. The jump count is stochastically dominated by a Poisson process of
rate \(r_+\), which proves nonexplosion, and invariance gives
\(r_-\leq\omega=\pi(r)\leq r_+\).

It remains to pass from the gap at \(s=0\) to the locally uniform tilted
decomposition. In the second case, define
\[
 (\mathcal Jf)(v)=\int_E W(v'|v)f(v')\dd v'.
\]
Then \(\|\mathcal J\|\leq r_+\). In the finite-state case the same
definition, with a sum in place of the integral, is a bounded matrix. Thus,
in either case,
\[
 \cL_s=\cL_0+(e^s-1)\mathcal J
\]
is an entire bounded operator family. Analytic perturbation of the isolated
simple eigenvalue gives an analytic eigenvalue \(\alpha(s)\) and an analytic
rank one Riesz projection \(\Pi_s\) \citep{Kato1995}. After reducing the
complex disc, the complementary spectrum of
\(\cL_s\) is separated from \(\alpha(s)\) by a uniform vertical gap. More
precisely, with \(Q_s=I-\Pi_s\), resolvent contours \(\Gamma_s\) can be
chosen with uniformly bounded lengths and resolvents so that they surround
\(\sigma(\cL_s|_{Q_s\mathcal X_{\mathbb C}})\) and satisfy
\[
 \sup_{z\in\Gamma_s}\Re z\leq\Re\alpha(s)-\gamma
\]
for some \(\gamma>0\). The Dunford integral then gives
\[
\begin{split}
 e^{t\cL_s}Q_s
 &=\frac1{2\pi i}\int_{\Gamma_s}
 e^{tz}(z-\cL_s)^{-1}Q_s\dd z,\\
 \|e^{t\cL_s}Q_s\|
 &\leq Ce^{t\{\Re\alpha(s)-\gamma\}},
\end{split}
\]
uniformly on the smaller disc. This proves
\[
 e^{t\cL_s}=e^{t\alpha(s)}\Pi_s+\mathcal R_s(t)
\]
with the asserted remainder bound, in both cases of the proposition.

Finally, \(\Pi_s\bfone\) remains nonzero near \(s=0\). With the invariant
functional \(\pi\) fixed, set
\[
 r_s=\Pi_s\bfone,\qquad
 \ell_s(f)=\frac{\pi(\Pi_sf)}{\pi(r_s)}.
\]
These functions are analytic, \(\Pi_s=r_s\otimes\ell_s\), and
\(\ell_s(r_s)=\ell_s(\bfone)=1\). At \(s=0\), \(r_0=\bfone\) and
\(\ell_0=\pi\). Moreover,
\(\|r_s-\bfone\|_\infty\to0\), so for every initial probability law \(\nu\),
\[
 a(s)=\langle\nu,\Pi_s\bfone\rangle=\nu(r_s)
\]
is bounded away from zero after one final, \(\nu\)-independent reduction of
the disc.
\end{proof}

In either setting of Proposition~\ref{prop:kinetic-spectral}, put
\(u_s(t,v)=\bE_v[e^{s\mathcal N(t)}]\). During the first \(\dd t\) units of
time the process either stays at \(v\) or jumps to a new state. Conditioning
on these possibilities gives
\[
\begin{split}
 u_s(t+\dd t,v)
 &=\{1-r(v)\dd t\}u_s(t,v)\\
 &\quad+\dd t\int_EW(v'|v)e^su_s(t,v')\dd v'\\
 &\quad+o(\dd t).
\end{split}
\]
Thus \(\partial_tu_s=\cL_su_s\) and \(u_s(0,\cdot)=\bfone\), so the
bounded generator uniqueness theorem gives
\(u_s(t,\cdot)=e^{t\cL_s}\bfone\). Integration against the initial law is
the jump process Feynman--Kac formula,
\[
 \bE_\nu[e^{s\mathcal N(t)}]
 =\langle\nu,e^{t\cL_s}\bfone\rangle.
\]
The assumed decomposition yields, locally uniformly in complex \(s\),
\[
\begin{gathered}
 \bE_\nu[e^{s\mathcal N(t)}]
 =e^{t\alpha(s)}a(s)\{1+\epsilon_t(s)\},\\
 \sup_{|s|<\varepsilon'}|\epsilon_t(s)|=O(e^{-\gamma t}).
\end{gathered}
\]
on every smaller disc on which \(a\) stays away from zero. For all
sufficiently large \(t\), the right side has no zero there, and the
analytic logarithm normalized at \(s=0\) is well-defined. It follows that
\[
 \frac1t\log\bE_\nu[e^{s\mathcal N(t)}]
 =\alpha(s)+\frac{\log a(s)+\log\{1+\epsilon_t(s)\}}t
\]
converges locally uniformly to \(\alpha(s)\). Consequently, the scaled
cumulant generating function,
\[
 \lim_{t\to\infty}\frac1t
 \log\bE_\nu[e^{s\mathcal N(t)}]
\]
equals \(\alpha(s)\). Cauchy's integral formula permits differentiation of
the locally uniform analytic limit and gives
\[
 \lim_{t\to\infty}\frac1t
 \cump_n\{\mathcal N(t)\}=\alpha^{(n)}(0).
\]
The Feynman--Kac and tilted generator principle is reviewed in
\citep{Touchette2009}; the semigroup and perturbation method is given in
\citep{EthierKurtz1986,Kato1995}. This is the jump process representation used
in \citep{ViscoVanWijlandTrizac2008}. Since each transition adds one collision,
the multiplier \(e^s\) appears only in the gain term. The source of memory is
the velocity left by the preceding collision, since that velocity determines
the next hazard.

For a tagged particle in a dilute equilibrium gas of hard-spheres, the
transition rates are fixed directly by the scattering rule. For a homogeneous
background with number density \(n_{\rm b}\), the Maxwellian density is
\[
 M_T(u)=\left(\frac{m}{2\pi k_{\rm B}T}\right)^{d/2}
 \exp\left(-\frac{m|u|^2}{2k_{\rm B}T}\right).
\]
For \(g=v-u\) and
\(\widehat\sigma\in\mathbb S^{d-1}\), set
\[
 v^\star=v-(g\cdot\widehat\sigma)\widehat\sigma .
\]
This is the tagged particle component of the equal mass collision rule
\eqref{eq:collision-rule}. With standard surface measure on
\(\mathbb S^{d-1}\), the tilted linear Boltzmann collision operator is
\begin{equation}\label{eq:linear-boltzmann-tilt}
\begin{split}
 &(\cL_s^{\rm B}f)(v)\\
 &\quad=n_{\rm b}a^{d-1}\int_{\bR^d}M_T(u)
 \int_{\mathbb S^{d-1}}(g\cdot\widehat\sigma)_+\\
 &\qquad\qquad\times\{e^sf(v^\star)-f(v)\}\dd\widehat\sigma\dd u.
\end{split}
\end{equation}
Here, \(x_+:=\max\{x,0\}\). In
Eq.~\eqref{eq:linear-boltzmann-tilt}, the scattering rule fixes the
transitions and their rates in the backward generator
\eqref{eq:tilted-generator}. The corresponding linear Boltzmann
equation has been derived from tagged particle dynamics of hard-spheres in the
Boltzmann--Grad scaling for the Rayleigh gas
\citep{MatthiesStoneTheil2018}. Visco, van Wijland and Trizac
\cite{ViscoVanWijlandTrizac2008} use this collision count tilt, with
the opposite sign convention for the transform parameter. Their first Sonine
projection of the eigenvalue problem gives
\eqref{eq:sonine1}--\eqref{eq:sonine3}, and hence
\(c_3-c_2<0\) in \eqref{eq:sonine-gap} within that kinetic closure.

Under the spatially homogeneous equilibrium Enskog closure
\citep{ChapmanCowling1970}, the contact pair correlation is a constant
\(g_{\rm eq}(a^+)>0\) and the tilted operator is
\begin{equation}\label{eq:homogeneous-enskog-tilt}
 \cL_s^{\rm E}=g_{\rm eq}(a^+)\cL_s^{\rm B}.
\end{equation}
In three dimensions, \(g_{\rm eq}(a^+)\) may be approximated by
\(g_{\rm CS}(\varphi)\) in \eqref{eq:carnahan-starling-contact}. Whenever the
dominant eigenvalue is defined as above,
\[
 \alpha_{\rm E}(s)=g_{\rm eq}(a^+)\alpha_{\rm B}(s),\qquad
 c_j^{\rm E}=g_{\rm eq}(a^+)c_j^{\rm B}.
\]
Multiplication by the positive contact factor rescales time and preserves
the normalized cumulant ratios and the sign of \(c_3-c_2\) within the
closure. The hard-sphere hazard in \eqref{eq:linear-boltzmann-tilt} grows
without bound as \(|v|\to\infty\). Applying the spectral argument to this
operator therefore requires estimates beyond the bounded rate condition of
Proposition \ref{prop:kinetic-spectral}. Its full spectral hypotheses and the
exact sign of \(c_3-c_2\) remain to be established here.

The next example isolates the effect of state-dependent waiting times in a
finite-state model. It gives a nonnegative pure jump count that reverses the
positive cluster inequality.

\begin{Prop}[Minimal two-state jump process counterexample]
\label{prop:two-state-latent}
Let \(E=\{1,2\}\), with transition rates
\[
 W(2|1)=a,\qquad W(1|2)=b,\qquad a,b>0,
\]
and let all other transition rates vanish. Let \(\mathcal N(t)\) denote the
number of transitions in \((0,t]\), and put
\[
 \chi=\frac{4ab}{(a+b)^2}.
\]
Then \(0<\chi\leq1\). For every initial probability law \(\nu\), the limits
\[
 c_j:=\lim_{t\to\infty}\frac{1}{t}
 \cump_j\{\mathcal N(t)\},\qquad j=1,2,3,
\]
exist and satisfy
\begin{align}
 c_1={}&\frac{a+b}{2}\chi,\nonumber\\
 c_2={}&(a+b)\chi\left(1-\frac{\chi}{2}\right),\nonumber\\
 c_3={}&(a+b)\chi
 \left(2-3\chi+\frac{3\chi^2}{2}\right).
\label{eq:two-state-cumulants}
\end{align}
In particular,
\begin{equation}\label{eq:two-state-gap}
\begin{gathered}
 c_3-c_2=\frac{(a+b)\chi}{2}(1-\chi)(2-3\chi)<0\\
 \text{if}\qquad \frac23<\chi<1.
\end{gathered}
\end{equation}
For example, \(a=2b\) gives
\[
 \frac{c_2}{c_1}=\frac{10}{9},
 \qquad
 \frac{c_3}{c_1}=\frac{28}{27}.
\]
\end{Prop}

\begin{proof}
The arithmetic--geometric mean inequality gives \(0<\chi\leq1\). Since the
chain is finite and irreducible, Proposition \ref{prop:kinetic-spectral}
applies. With the states ordered as \(1,2\), the tilted generator is
\[
 \cL_s=
 \begin{pmatrix}
 -a&ae^s\\
 be^s&-b
 \end{pmatrix}.
\]
The eigenvalue that equals zero at \(s=0\) is
\[
\begin{split}
 \alpha(s)&=\frac{-(a+b)+\sqrt{(a-b)^2+4abe^{2s}}}{2}\\
 &=\frac{a+b}{2}\left[\sqrt{1+\chi(e^{2s}-1)}-1\right].
\end{split}
\]
Its expansion at zero is
\[
\begin{split}
 \alpha(s)=\frac{a+b}{2}\Bigg[&\chi s+
 \left(\chi-\frac{\chi^2}{2}\right)s^2\\
 &+\left(\frac{2\chi}{3}-\chi^2+\frac{\chi^3}{2}\right)s^3
 +O(s^4)\Bigg].
\end{split}
\]
The tilted generator formula gives \(c_j=\alpha^{(j)}(0)\). Differentiating
the last display proves \eqref{eq:two-state-cumulants}. Subtraction gives
\eqref{eq:two-state-gap}, and the stated example follows by taking
\(\chi=8/9\).
\end{proof}

Every event in Proposition \ref{prop:two-state-latent} switches the state,
so the holding rates alternate. The two-states serve as abstract labels. A
discretization of \eqref{eq:linear-boltzmann-tilt} would also have to count
collisions whose initial and final velocities lie in the same cell. The
example establishes that state-dependent holding rates can produce
\(c_3<c_2\) for a nonnegative count. The positive Hawkes cluster construction
therefore imposes a restriction on the dependence between events.

\subsection{A clipped signed pair process}

Let \(h_{\rm same},h_{\rm share},h_{\rm dis}\) now be signed integrable
kernels, and construct \(h_{ef}\) from the three intersection relations as in
\eqref{eq:pair-hawkes}. A nonlinear extension is
\begin{equation}\label{eq:clipped-model}
 \lambda_e(t)=\left[\mu+\sum_{f\in\Etwo}
 \int_{(-\infty,t)}h_{ef}(t-s)N_f(\dd s)\right]_+.
\end{equation}
A negative part of a kernel reduces the intensity after an event. Clipping
at zero keeps the rate nonnegative. This can represent a reduction in
activity at short lags, such as the interval immediately after two hard
spheres separate. The next proposition gives a sufficient stationarity
condition.

\begin{Prop}[A sufficient stationarity condition]\label{prop:clipped}
Define the nonnegative matrix \(A\) as follows.
\[
 A_{ef}=\int_0^\infty|h_{ef}(t)|\dd t.
\]
If \(\rho(A)<1\), then \eqref{eq:clipped-model} has a unique stationary,
nonexplosive solution with finite mean intensity. If the three kernels and
the baseline are exchangeable, its stationary law is exchangeable. Moreover,
the integrated Lipschitz matrix of the nonlinear intensity map has norm
strictly less than one in a weighted coordinate norm. In the exchangeable
case, put
\[
 a_\alpha=\int_0^\infty|h_\alpha(t)|\dd t,
 \qquad \alpha\in\{{\rm same},{\rm share},{\rm dis}\}.
\]
Then the condition \(\rho(A)<1\) is equivalent to the explicit inequality
\[
 a_{\rm same}+2(N-2)a_{\rm share}
 +\binom{N-2}{2}a_{\rm dis}<1.
\]
\end{Prop}

\begin{proof}
In the exchangeable case, every row of the nonnegative matrix \(A\) has sum
\[
 a_{\rm abs}:=a_{\rm same}+2(N-2)a_{\rm share}
 +\binom{N-2}{2}a_{\rm dis}.
\]
Hence \(A\bfone_M=a_{\rm abs}\bfone_M\), while
\(\rho(A)\leq\|A\|_\infty=a_{\rm abs}\); consequently
\(\rho(A)=a_{\rm abs}\). In this case the ordinary coordinate supremum norm
already gives the claimed Lipschitz bound whenever \(a_{\rm abs}<1\).

For a general nonnegative matrix \(A\), the contraction can be made explicit.
The link \(x\mapsto[x]_+\) is nonnegative and Lipschitz with constant one.
Choose \(q\) with \(\rho(A)<q<1\) and set
\[
 w=(I_M-A/q)^{-1}\bfone_M
   =\sum_{n\geq0}(A/q)^n\bfone_M.
\]
Then \(w\) has strictly positive coordinates. Moreover
\[
 Aw=q(w-\bfone_M)\leq qw.
\]
Thus, in the weighted coordinate norm
\(\|x\|_w=\max_e|x_e|/w_e\), the integrated Lipschitz matrix of the
intensity map satisfies
\[
 \|A\|_w
 =\max_e\frac{\sum_fA_{ef}w_f}{w_e}
 =\max_e\frac{(Aw)_e}{w_e}
 \leq q<1.
\]
Indeed, for any two input vectors \(x,y\),
\[
 |[\mu+x_e]_+-[\mu+y_e]_+|\leq|x_e-y_e|,
\]
and integration of the coordinatewise kernel bounds produces precisely the
matrix \(A\). The contraction criterion for
nonlinear multivariate point processes in Theorem 7 of
\citep{BremaudMassoulie1996} gives existence, uniqueness, finite mean
intensity and stability. Finite mean intensity implies
\(\bE N_e(0,t]<\infty\) for every bounded interval, hence local finiteness
and nonexplosion. If the coordinates are permuted by a particle
permutation, the transformed process satisfies the same stochastic intensity
equation. Uniqueness in law therefore implies exchangeability.
\end{proof}

Clipping individual coordinates also changes how perturbations pass between
the three modes.

\begin{Rem}[Modal coupling under the clipped link]
Let \(Z(t)\) be the vector inside the brackets in
\eqref{eq:clipped-model}, and put \(\Phi(x)=([x_e]_+)_{e\in\Etwo}\). At any
\(x\in\bR^M\) with \(x_e\neq0\) for every \(e\), the derivative of the link is
\[
 D\Phi(x)=D_x:=\operatorname{diag}
 \{\mathbf 1_{\{x_e>0\}}:e\in\Etwo\}.
\]
For every lag \(u\), the signed kernel matrix belongs to the Johnson algebra
and can therefore be written as
\[
 H(u)=h_G(u)P_G+h_P(u)P_P+h_C(u)P_C,
\]
where \(h_G(u),h_P(u),h_C(u)\) are its three modal eigenvalues. Hence, for
\(a,b\in\{G,P,C\}\),
\begin{equation}\label{eq:clipped-modal-transfer}
 P_aD\Phi(x)H(u)P_b=h_b(u)P_aD_xP_b.
\end{equation}
If no coordinate is clipped, then \(D_x=I_M\), and the off diagonal blocks in
\eqref{eq:clipped-modal-transfer} vanish. Suppose instead that exactly one
coordinate \(e\) is clipped, and let \(\mathsf e_e\) be its coordinate vector.
Then \(D_x=I_M-\mathsf e_e\mathsf e_e^\top\), and, for \(a\neq b\),
\[
 P_aD_xP_b=-(P_a\mathsf e_e)(P_b\mathsf e_e)^\top\neq0.
\]
Indeed, the squared norms of
\(P_G\mathsf e_e,P_P\mathsf e_e,P_C\mathsf e_e\) are the diagonal entries of
the corresponding projectors, namely
\[
 \frac1M,\qquad \frac2N,\qquad \frac{N-3}{N-1},
\]
which are positive for \(N\geq4\). For this one-coordinate clipping mask,
the transfer block from mode \(b\) to each other mode is nonzero whenever
\(h_b(u)\neq0\).

Under an exchangeable stationary law, the unconditional cross-mode
covariances vanish whenever the second moments are finite. For
\(\Delta>0\) such that \(\bE\|N(0,\Delta]\|^2<\infty\), define
\[
 \Sigma_\Delta=\operatorname{Cov}\{N(0,\Delta]\}.
\]
Exchangeability implies that the entries of \(\Sigma_\Delta\) depend only on
whether two edges are equal, share one endpoint or are disjoint. Therefore
\(\Sigma_\Delta\) belongs to the Johnson algebra and has the exact form,
\begin{equation}\label{eq:clipped-window-covariance}
\begin{split}
 \Sigma_\Delta
 &=v_G(\Delta)P_G+v_P(\Delta)P_P+v_C(\Delta)P_C,\\
 v_a(\Delta)&=\frac{\operatorname{tr}(P_a\Sigma_\Delta)}{\dim U_a}.
\end{split}
\end{equation}
Each \(v_a(\Delta)\) may depend on all three signed kernels through the
stationary clipping mask. Exchangeability fixes the three covariance
subspaces. Determining their coefficients requires the law of the nonlinear
process, in place of the scalar calculation used in
\eqref{eq:finite-window-factor}.
\end{Rem}

With finite second moments, Eq.~\eqref{eq:clipped-window-covariance}
follows from exchangeability.
The coefficients in \eqref{eq:hawkes-covariance} and the proof of Theorem
\ref{thm:cumulant} use the positive linear Hawkes structure. Renewal and
likelihood questions for inhibitory Hawkes processes are studied in
\citep{CostaGrahamMarsalleTran2020,BonnetHerreraSangnier2021}.

\subsection{Other kinetic descriptions with memory}

Memory also appears in other kinetic descriptions. Projection methods lead
to generalized kinetic equations with memory kernels
\citep{Zwanzig1961,Mori1965}. Moffat \cite{Moffat1974} gave a
stochastic interpretation of model kinetic equations with memory in which
collisions form a renewal process, while equations of Enskog type with memory
have been developed for granular gases
\citep{BorovchenkovaGerasimenko2014}. Continuous time random walks give
another description based on waiting times
\citep{MontrollWeiss1965,MetzlerKlafter2000}. These approaches place memory in
a waiting time law or in an evolution operator, and require a separate
cumulant analysis from Theorem \ref{thm:cumulant}. For the latent velocity
operator \eqref{eq:tilted-generator}, retaining the velocity gives a Markov
description. The collision count obtained by omitting that velocity can
retain dependence on its past.

\section{Conclusion}\label{sec:conclusion}

We have proved that every nonnegative subset count of a stationary
subcritical positive Hawkes process with finitely many types has long time
cumulant rates that are nondecreasing in their order. The reason is the
nonnegative integer contribution of each independent immigrant cluster.
The proof includes clusters crossing either endpoint of the observation
window, and applies to total and tagged particle counts.

The first Sonine approximation for hard particles gives a third cumulant
rate smaller than the second in every dimension at least two. This ordering
rules out an exact match to that kinetic benchmark by a positive Hawkes
subset count. Establishing the same sign for exact hard-sphere dynamics
remains a separate problem. The planar simulations support the reversed
ordering at the longer observation window. The three-dimensional
calculation leaves that ordering unresolved.

For pair counts, the Johnson decomposition separates total activity,
differences between tagged counts and cycle contrasts. It gives the
stationarity condition, mean-rate calibration and covariance formulas.
With exponential kernels, the variance can also be computed for each finite
observation window. The simulations agree with the exact finite-window
covariance factors across the displayed particle numbers and window
widths. Cycle contrasts use signed weights, so their cumulants
require a separate analysis from the subset count inequality.

The two-state example shows how the state left by one event can change the
next waiting time and produce \(c_3<c_2\). In the kinetic description, the
first Sonine closure of the tilted linear Boltzmann operator has this
ordering. Its full spectral hypotheses and exact cumulant sign remain to
be established here. For the clipped signed process, we obtain a sufficient
stationarity condition; determining its cumulants remains a further problem.

\begin{acknowledgments}
I want to thank my professor and longtime friend Marcos Ley Koo for his advice
throughout the years and for the nice conversations that inspired this
research. I also want to thank my parents, Guadalupe and Mart\'in, my brother
Omar, and my partner Kayoko for their unwavering support in all my pursuits.

This research received no external funding.

OpenAI Codex (GPT-6) assisted with the organization and language revision of
this manuscript under the author's direction. The preparation included
checks of mathematical statements, references and the compiled document.
The author is responsible for the content.
\end{acknowledgments}

\section*{Data availability}
The data and reproducibility materials supporting this study are openly
available in Zenodo at \url{https://doi.org/10.5281/zenodo.21837445}
\cite{HernandezRuiz2026Data}.

\bibliography{references}
\end{document}